\documentclass[a4paper]{amsart}

\usepackage{amsmath,amssymb}
\usepackage{xcolor}
\usepackage{pictexwd}
\usepackage[colorlinks=true,allcolors=blue]{hyperref}

\newtheorem{theorem}{Theorem}[section]
\newtheorem{proposition}[theorem]{Proposition}
\newtheorem{corollary}[theorem]{Corollary}
\newtheorem{lemma}[theorem]{Lemma}
\newtheorem{algorithm}[theorem]{Algorithm}

\newtheorem{preremark}[theorem]{Remark}
\newtheorem{predefinition}[theorem]{Definition}
\newtheorem{preexample}[theorem]{Example}
\newtheorem{prenotation}[theorem]{Notation}
\newtheorem{preconjecture}[theorem]{Conjecture}

\newenvironment{remark}{\begin{preremark}\rm}{\end{preremark}}
\newenvironment{definition}{\begin{predefinition}\rm}
{\end{predefinition}}
\newenvironment{example}{\begin{preexample}\rm}{\end{preexample}}
\newenvironment{conjecture}{\begin{preconjecture}\rm}
{\end{preconjecture}}

\allowdisplaybreaks
\def\OO{{\mathcal{O}}}

\def\AA{\mathbb{A}}
\def\NN{\mathbb{N}}
\def\ZZ{\mathbb{Z}}
\def\QQ{\mathbb{Q}}

\newcommand{\M}{{\mathfrak{M}}}

\newcommand{\m}{{\mathfrak{m}}}

\let\epsilon=\varepsilon

\def\phi{{\varphi}}
\let\Psi=\varPsi
\let\Phi=\varPhi
\let\theta=\vartheta
\let\rho=\varrho

\def\LT{\mathop{\rm LT}\nolimits}

\def\ND{\mathop{\rm ND}\nolimits}
\def\AR{\mathop{\rm AR}\nolimits}
\def\HF{\mathop{\rm HF}\nolimits}

\def\Mat{\mathop{\rm Mat}\nolimits}

\def\Supp{\mathop{\rm Supp}\nolimits}
\def\Spec{\mathop{\rm Spec}\nolimits}
\def\Hilb{\mathop{\rm Hilb}\nolimits}
\def\Cot{\mathop{\rm Cot}\nolimits}
\def\Syz{\mathop{\rm Syz}\nolimits}

\def\wt{\mathop{\rm wt}\nolimits}

\newcommand{\lin}{\mathop{\rm lin}\nolimits}
\newcommand{\Lin}{\mathop{\rm Lin}\nolimits}
\newcommand{\edim}{\mathop{\rm edim}\nolimits}
\newcommand{\BO}{\mathbb{B}_{\mathcal{O}}}
\newcommand{\BOhom}{\mathbb{B}_{\mathcal{O}}^{\rm hom}}
\newcommand{\BOdf}{\mathbb{B}_{\mathcal{O}}^{\rm df}}

\newcommand{\Cint}{C^{\rm int}}
\newcommand{\Crim}{C^{\rm rim}}
\newcommand{\Cexp}{C^{\rm exp}}

\newcommand{\Cplus}{C_+}
\newcommand{\Cnull}{C_0}

\newcommand{\OOrim}{\OO^{\rm rim}}
\newcommand{\OOint}{\OO^{\rm int}}

\def\TTTo#1\mathop{\xrightarrow{\hspace*{1cm}}^{^{\mkern-70mu (#1}}}

\def\tfrac #1#2{{\textstyle\frac{#1}{#2}}}
\def\tsum_#1^#2{{\textstyle\sum\limits_{#1}^{#2}}}
\def\tbinom #1#2{{\textstyle\binom{#1}{#2}}}

\newcommand{\gray}[1]{\textcolor{gray}{#1}}

\def\cocoa{\mbox{\rm
  C\kern-.13em o\kern-.07 em C\kern-.13em o\kern-.15em A}}
\def\apcocoa{\mbox{\rm
A\kern-0.13em p\kern -0.07em C\kern-.13em o\kern-.07 em C\kern-.13em
o\kern-.15em A}}

\begin{document}

\title{Re-Embeddings of Special Border Basis Schemes}

\author{Martin Kreuzer}
\address{Fakult\"at f\"ur Informatik und Mathematik, 
Universit\"at Passau, D-94030 Passau, Germany}
\email{martin.kreuzer@uni-passau.de}

\author{Lorenzo Robbiano}
\address{Dipartimento di Matematica, Universit\`a di Genova,
Via Dodecaneso 35,
I-16146 Genova, Italy}
\email{lorobbiano@gmail.com}

\dedicatory{This paper is dedicated to the memory of Wolmer Vasconcelos.}


\begin{abstract}
Border basis schemes are open subschemes of the Hilbert scheme of $\mu$ points in 
an affine space $\mathbb{A}^n$. They have easily describable systems of generators of their vanishing ideals
for a natural embedding into a large affine space $\mathbb{A}^{\mu\nu}$. Here we bring together
several techniques for re-embedding affine schemes into lower dimensional spaces which
we developed in the last years. We study their efficacy for some special
types of border basis schemes such as MaxDeg border basis schemes, L-shape and simplicial border basis schemes,
as well as planar border basis schemes. A particular care is taken to make these re-embeddings
efficiently computable and to check when we actually get an isomorphism with $\mathbb{A}^{n\mu}$, i.e.,
when the border basis scheme is an affine cell.
\end{abstract}






\maketitle

\tableofcontents

%
%

\section{Introduction}
\label{sec1}

The second author enjoys remembering when, using his supreme irony, 
his deep and profound voice, Wolmer Vasconcelos liked to declaim the 
following last stanza of the poem ``Sat\'elite'' by Manuel Bandeira.

\begin{flushright}
\begin{tabular}{ll}
{\it Fatigado de mais-valia,} & Tired of the added value, \\
{\it Gosto de ti assim:}      & I like you this way:  \\
{\it Coisa em si,}            & Thing in itself,  \\
{\it -- Sat\'elite}           & -- Satellite
\end{tabular}
\end{flushright}

It is not clear how much added value is provided by the theory of border basis
schemes. Certainly they form open subsets covering the Hilbert scheme of points,
a subject which was initiated by A.\ Grothendieck in~\cite{Gro} and 
has been the topic of hundreds of papers already. In contrast to
Hilbert schemes, border basis schemes are defined by easily described polynomials, making them
eminently suitable for computer calculations. The authors introduced the term {\it border basis scheme}
when they studied flat deformations of border bases in their 2008 paper~\cite{KR3}.
However, viewed as open subschemes of $\Hilb^\mu(\mathbb{A}^n)$, they had been
examined previously by M.\ Haiman for the case $n=2$ in~\cite{Hai}, and subsequently by 
M.\ Huibregtse in the general case (see~\cite{Hui1,Hui2,Hui3}).
Notwithstanding the fact that they yield an open covering of the Hilbert scheme of points, 
after working on border basis schemes for about 20 years, the authors view them as beautiful 
objects in themselves, even without an added value.

For schemes in algebraic geometry, it is always an important task to determine into which
affine or projective spaces they can be embedded. For instance, a classical theorem states that a
non-singular projective curve can always be re-embedded into~$\mathbb{P}^3$. Thus it is a natural
question to ask for re-embeddings of border basis schemes into low dimensional affine spaces.
Another interesting problem is to decide whether a border basis scheme is actually isomorphic
to some affine space. In this case, being an open subscheme of the Hilbert scheme, 
it is called an {\it affine cell}.

When the authors began to study the topic of re-embedding border basis schemes into
lower dimensional affine spaces, they soon ran into apparently insurmountable computational
problems. The huge number of indeterminates involved in their defining equations
made the classical calculation of elimination ideals using Gr\"obner bases all but impossible.
At this point the beautiful structural properties of the defining equations of a border basis 
scheme came to the rescue. Many of them are of the form $x_i - f(x_1, \dots, \widehat{x_i},
\dots, x_n)$ which readily allows us to eliminate~$x_i$. Naturally, the situation gets
much less obvious and, indeed, quite tricky, if we want to eliminate several indeterminates
simultaneously in this way. Thus the method of {\it elimination 
by substitution} was born. Always having the ideals defining border basis schemes in mind,
we developed this method in a series of papers (see~\cite{KLR1,KLR2,KLR3,KR5,AKL}), 
some of them together with L.N.~Long.

Finally, after this detour, we come back to the original task of
re-embedding border basis schemes and try to apply the method of
elimination by substitution, or $Z$-separating re-embeddings, as we originally called it,
to the fullest. Thus, even if many results in this paper look like easy applications of
the theory of $Z$-separating re-embeddings, don't let that fool you: in themselves, they
are not obvious. 

Let us  briefly go through the main dishes we serve up here. Given an order ideal
of terms~$\OO$, the border basis scheme $\BO$ parametrizes all 0-dimensional ideals~$I$ in
the polynomial ring $P=K[x_1,\dots,x_n]$ over a field~$K$ such that~$\OO$ represents a
$K$-basis of~$P/I$.

In Section~\ref{sec2} we start by recalling two ways to generate the vanishing ideal $I(\BO)$
in $K[C] = K[c_{11},\dots,c_{\mu\nu}]$ for the natural embedding of~$\BO$.
The arrow grading is introduced, and the indeterminates~$C$ are clustered in several ways: 
rim and interior indeterminates, as well as
exposed and non-exposed indeterminates. The exposed indeterminates are characterized
in different ways which explain their name (see Proposition~\ref{prop:CharExposed}).

Section~\ref{sec3} first recapitulates the method of $Z$-separating re-embeddings and then
discusses MaxDeg border basis schemes. They correspond to MaxDeg order ideals, i.e., to order ideals
which have the property that no term in~$\OO$ has a larger degree than any term in the 
border~$\partial\OO$ of~$\OO$. MaxDeg order ideals are characterized in several ways
(see Proposition~\ref{prop-CharMaxDeg}) and it is shown that MaxDeg border basis schemes satisfy 
the assumptions at the core of~\cite{KR5}: the total arrow grading is non-negative
and their vanishing ideal $I(\BO)$ intersects $K[C]_0$ trivially 
(see Proposition~\ref{prop-MaxDegBasis}). Thus we can find optimal $Z$-separating
re-embeddings effectively and detect when~$\BO$ is an affine cell (see Proposition~\ref{prop-MaxDegReEmbed}).

A particular MaxDeg border basis scheme has been somewhat of a mystery for quite some time, namely
the L-shape border basis scheme which we look at in Section~\ref{sec4}. It corresponds to the order ideal
$\OO=\{1,y,x,y^2,x^2\}$ in two indeterminates, has dimension~10, and is naturally embedded 
in~$\mathbb{A}^{25}$. Using the best $Z$-separating re-embedding, we can realize it as a
codimension~2 smooth complete intersection in~$\mathbb{A}^{12}$ (see Proposition~\ref{prop-L-Shape-Zsep}). 
For a long while, it was believed that this is all one can do. However, using a new
technique based on the Unimodular Matrix Problem developed in~\cite[Sect.~6]{KR5}, we are
now able to show that $\BO\cong\mathbb{A}^{10}$, i.e., that the L-shape border basis scheme 
is an affine cell (see Proposition~\ref{prop-theLshape}).

Another special case, namely simplicial order ideals, i.e., order ideals of the form 
$\OO = \{ t\in\mathbb{T}^n \mid \deg(t)\le d\}$ for some $d\ge 1$, is treated in Section~\ref{sec5}.
Since they have a positive total arrow grading (see Corollary~\ref{cor-CharSimplicial}), 
we know already that we can get an {\it optimal} $Z$-separating re-embedding, i.e., a re-embedding 
which reaches the embedding dimension. Here we go one step further: we construct a concrete
re-embedding based on a $Z$-separating tuple of polynomials consisting of some
of the natural generators of~$I(\BO)$ (see Proposition~\ref{prop-simplicialZ}). Furthermore, 
we deduce that simplicial border basis schemes are affine cells for $n=2$ and singular at
the monomial point for $n>2$ (see Proposition~\ref{prop-AffineOrSing}).

In Sections~\ref{sec6} and~\ref{sec7}, we examine another special kind of border basis schemes,
namely planar border basis schemes. Like the L-shape border basis scheme, they correspond to order 
ideals in $n=2$ indeterminates
and are known to be smooth and irreducible schemes of dimension $2\mu$, where $\mu=\#\OO$.
Our first major contribution in Section~\ref{sec6} is the Weight Assignment Algorithm~\ref{alg-WeightAssign}.
It allows us to define weights for the indeterminates and choose natural generators
of~$I(\BO)$ such that the corresponding $Z$-separating re-embedding eliminates
all non-exposed indeterminates~$c_{ij}$. Then, in Theorem~\ref{thm-elimNex} 
and Corollary~\ref{cor-exposedgenerate}, we reap the rewards: we eliminate all 
non-exposed indeterminates and show that the exposed indeterminates generate the coordinate ring
of a planar border basis scheme. We also compare this approach to the methods described
in~\cite[Sect.~6 and~7]{Hui1}, and we provide a variant of~\cite[Alg.~5.7]{KLR3} which computes
the optimal $Z$-separating re-embeddings in the planar case if they exist (see Algorithm~\ref{alg-OptPlanarZ}).

The final Section~\ref{sec7} deals with two special types of planar border basis schemes.
The first type are planar box border basis schemes defined by order ideals of the form 
$\OO = \{ x^i y^j \mid 0\le i< a,\; 0\le j< b\}$ with $a,b>0$. For them, it suffices to
eliminate the non-exposed indeterminates $c_{ij}$ in order to get affine cells 
(see Proposition~\ref{prop-PlanarBox}).
The second type are planar MaxDeg border basis schemes. For them, we know already that they
are affine cells (at least if $K$ is perfect), but the isomorphism $\BO \cong \mathbb{A}^{2\mu}$
may involve some complicated component derived from the Unimodular Matrix Problem (see Proposition~\ref{prop-theLshape}).
Thus it is a natural question to ask when a $Z$-separating re-embedding suffices to construct such
an isomorphism. In Conjecture~\ref{conj-newconj}, we suggest an answer to this question which depends
only on the number of segments of maximal degree of~$\OO$ (as introduced in Definition~\ref{def-OdSegment}).
This conjecture is supported by extensive numerical evidence and constitutes a fitting conclusion to this
{\it tour de force} through the world of re-embeddings of special border basis schemes.

The definitions and results in this paper are illustrated by suitable examples
These examples were computed using the computer algebra system \cocoa, 
to which several new or updated packages were added, for instance the packages
{\tt MaxDegBBS.cpkg5} and {\tt NonNegGrading.cpkg5}, as well as
{\tt QuillenSuslin.cpkg5}, where the latter package was implemented by L.N.\ Long.
The notation and terminology we adhere to are based on the books~\cite{KR1} and~\cite{KR2}.

\bigskip\bigbreak
%
%

\section{Generalities on Border Basis Schemes}
\label{sec2}

In this paper we study re-embeddings of border basis schemes.
Thus we begin by recalling the definition of border basis schemes
and some generalities about them. A border basis scheme is always given 
with respect to some order ideal of terms which, in turn, is a concept defined as follows.

Let $K$ be a field and $P=K[x_1, \dots, x_n]$. The set of terms in~$P$
is denoted by $\mathbb{T}^n = \{ x_1^{\alpha_1}\cdots x_n^{\alpha_n}
\mid \alpha_i \ge 0\}$. A finite set of terms $\OO = \{t_1,\dots,  t_\mu\}$ in $\mathbb{T}^n$
is called an {\bf order ideal}, if $t\in\OO$ and $t'\mid t$ implies $t'\in \OO$,
i.e., if every term dividing a term in~$\OO$ is again in~$\OO$.
The name ``order ideal'' is derived from the fact that such sets are poset ideals
in the poset $\mathbb{T}^n$ which is partially ordered by divisibility.

Given an order ideal $\OO = \{t_1,\dots,  t_\mu\}$, we call
$\partial\OO = (x_1\OO \cup \cdots \cup x_n\OO)
\setminus \OO$ the {\bf border} of~$\OO$.
We always write $\partial\OO=\{ b_1, \dots, b_\nu\}$ with $b_i\in\mathbb{T}^n$.

In this setting, the $\OO$-border basis scheme can be 
introduced as follows.

\begin{definition}
Let $\OO = \{t_1, \dots, t_\mu\}$ be an order ideal in~$\mathbb{T}^n$, and let
$\partial\OO = \{b_1, \dots, b_\nu\}$ be its border.
\begin{enumerate}
\item[(a)] For $i=1,\dots,\mu$ and $j=1,\dots,\nu$, let $c_{ij}$ be a new indeterminate.
Let $C=(c_{11}, c_{12}, \dots, c_{\mu\nu})$ be the tuple of all these indeterminates,
and let $K[C] = K[c_{11},\dots,c_{\mu\nu}]$.
Then the set of polynomials $G=\{ g_1,\dots,  g_\nu\}$ in $K[C][x_1,\dots,x_n]$, where 
$$
g_j \;=\; b_j - c_{1j}\, t_1 - \cdots - c_{\mu j} t_\mu
$$ 
for $j=1,\dots,\nu$, is called the {\bf generic $\OO$-border prebasis}.

\item[(b)] For $r=1,\dots,n$, the matrix $\mathcal{A}_r = (a_{ij}^{(r)}) \in \Mat_\mu(K[C])$
whose entries are given by
$$
a_{ij}^{(r)} \;=\; \begin{cases}
\delta_{im} & \hbox{\rm if\ } x_r t_j = t_m, \\
c_{im} & \hbox{\rm if\ } x_r t_j = b_m,
\end{cases}
$$
is called the $r$-th {\bf generic multiplication matrix} for~$\OO$.

\item[(c)]  Let $I(\BO) \subseteq K[C]$ be the ideal generated
by the entries of the commutator matrices $\mathcal{A}_i \mathcal{A}_j - \mathcal{A}_j
\mathcal{A}_i$, where $1\le i < j \le n$. Then the subscheme of $\mathbb{A}^{\mu \nu}$
defined by $I(\BO)$ is called the {\bf $\OO$-border basis scheme} and is denoted by~$\BO$.

\item[(d)] The coordinate ring of~$\BO$ is denoted by $B_\OO = K[C] / I(\BO)$.

\end{enumerate}
\end{definition}

The generators of the ideal $I(\BO)$ can also be constructed in the following way.

\begin{remark}\label{NaturalGens}
As in~\cite[Def. 6.4.33]{KR2} and~\cite[Def. 6.1]{KLR3}, we say that two border terms~$b_j$ and~$b_{j'}$
are {\bf next-door neighbors} if $b_j=x_k b_{j'}$ for some $k\in \{1,\dots,n\}$
and {\bf across-the rim neighbors} if $b_j= x_k t_m$ and $b_{j'} = x_\ell t_m$ for some
$k,\ell\in \{1,\dots,n\}$ and $m\in \{1,\dots,\mu\}$.

Following~\cite[Sect.~4]{KR3} or~\cite[Def.~6.1]{KLR3}, we construct the set
of {\bf next-door generators} $\ND_\OO$ as the union of the entries of the tuples $\ND(j,j')$,
where $b_j$ and $b_{j'}$ are next-door neighbors. Similarly, we construct
the set of {\bf across-the-rim generators} $\AR_\OO$ as the union of the entries 
of the tuples $\AR(j,j')$, where~$b_j$ and~$b_{j'}$ are across-the rim neighbors.

Finally, we take the union $\ND_\OO \cup \AR_\OO$ and get $I(\BO) = \langle \ND_\OO \cup \AR_\OO \rangle$.
This set is also called the set of {\bf natural generators} of~$I(\BO)$.
\end{remark}

The ideal $I(\BO)$ is generated by polynomials of degree two which are contained
in the maximal ideal $\M = \langle c_{ij} \mid i=1,\dots,  \mu;\; j=1,\dots,  \nu\rangle$
of~$K[C]$. The linear and quadratic parts of the natural generators of $I(\BO)$
can be described explicitly (see~\cite[Cor. 2.8]{KSL} and~\cite[Props. 6.3 and 6.6]{KLR3}).

A further important property of the ideals $I(\BO)$ is that they are homogeneous
with respect to the arrow grading and hence the total arrow grading 
(see~\cite[Lemma~3.4]{KSL} and~\cite[Prop.~6.3]{KLR3}).
Here the {\bf arrow grading} is the $\ZZ^n$-grading on~$K[C]$ 
defined by the matrix $A\in\Mat_{n,\mu\nu}(\ZZ)$
such that $\deg_A(c_{ij}) = \log(b_j) - \log(t_i)$.
(Recall that, for a term $t=x_1^{\alpha_1}
\cdots x_n^{\alpha_n}$, we let $\log(t)=(\alpha_1,\dots,  \alpha_n)$.)
Moreover, the {\bf total arrow grading} is the $\ZZ$-grading on~$K[C]$ 
defined by the tuple $W\in \ZZ^{\mu\nu}$ with $\deg_W(c_{ij}) = \deg(b_j) - \deg(t_i)$.

The subtuple $\Cplus = (c_{ij} \mid \deg_W(c_{ij})>0)$ of~$C$
is called the tuple of {\bf positive indeterminates}, and the subtuple
$\Cnull = (c_{ij} \mid \deg_W(c_{ij})=0 )$ is called the tuple of
{\bf zero degree indeterminates}. 

As for the arrow degrees of the indeterminates~$c_{ij}$, we have the
following observation (see also~\cite[Lemma 4.1.1]{Hui2}).

\begin{proposition}\label{prop-poscomp}
For all $i\in\{1,\dots,\mu\}$ and $j\in\{1,\dots,\nu\}$, the arrow degree
$\deg_A(c_{ij})$ has at least one positive component.
\end{proposition}

\begin{proof}
Suppose that $i\in\{1,\dots,\mu\}$ and $j\in\{1,\dots,\nu\}$
such that all components of $\log(b_j) -\log(t_i)$ are negative.
This means that $b_j \mid t_i$, and therefore the fact that~$\OO$
is an order ideal implies $b_j\in\OO$, a contradiction.
\end{proof}

For later usage, we subdivide the tuple of indeterminates $C=(c_{ij})$ as follows.

\begin{definition}\label{def-rim-and-interior}
Let $\OO = \{t_1, \dots, t_\mu\}$ be an order ideal in~$\mathbb{T}^n$, and let
$\partial\OO = \{b_1, \dots, b_\nu\}$ be its border.
\begin{enumerate}
\item[(a)] For $i\in \{1,\dots, \mu\}$, the term~$t_i$ is called a {\bf rim term}
in~$\OO$ if there exists an indeterminate~$x_k$ such that $x_k t_i \in \partial\OO$.
Otherwise, the term~$t_i$ is called an {\bf interior term} of~$\OO$.
The set of all rim terms is denoted by $\OOrim$, and the set of all interior terms
is denoted by $\OOint$.

\item[(b)] For $i\in \{1,\dots, \mu\}$ and $j\in \{1,\dots, \nu\}$, the 
indeterminate $c_{ij}$ is called a {\bf rim indeterminate} if~$t_i$ is a rim term
in~$\OO$. The set of all rim indeterminates will be denoted by~$\Crim$.

\item[(c)] For $i\in \{1,\dots, \mu\}$ and $j\in \{1,\dots, \nu\}$, the 
indeterminate $c_{ij}$ is called an {\bf interior indeterminate} if~$t_i$ is an
interior term of~$\OO$. The set of all interior indeterminates will be denoted
by~$\Cint$.

\end{enumerate}
\end{definition}

Let us look at these notions in a concrete case.

\begin{example}{\bf (The (2,2)-Box)}\label{ex-box22}\\
In the polynomial ring $P=K[x,y]$ over a field~$K$, 
we consider the order ideal $\OO = \{t_1,t_2,t_3,t_4\} = \{1,\, y,\, x,\, xy\}$ 
which we call the {\bf (2,2)-box}.
\begin{figure}[ht]
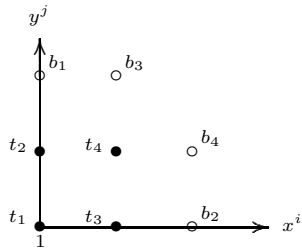

\centering{\makebox{
\beginpicture
   \setcoordinatesystem units <1cm,1cm>
   \setplotarea x from 0 to 3, y from 0 to 2.5
   \axis left /
   \axis bottom /
   \arrow <2mm> [.2,.67] from  2.5 0  to 3 0
   \arrow <2mm> [.2,.67] from  0 2  to 0 2.5
   \put {$\scriptstyle x^i$} [lt] <0.5mm,0.8mm> at 3.1 0.1
   \put {$\scriptstyle y^j$} [rb] <1.7mm,0.7mm> at 0 2.6
   \put {$\bullet$} at 0 0
   \put {$\bullet$} at 1 0
   \put {$\bullet$} at 0 1
   \put {$\bullet$} at 1 1

   \put {$\scriptstyle 1$} [lt] <-1mm,-1mm> at 0 0
   \put {$\scriptstyle t_1$} [rb] <-1.3mm,0.4mm> at 0 0
   \put {$\scriptstyle t_3$} [rb] <-1.3mm,0.4mm> at 1 0
   \put {$\scriptstyle t_2$} [rb] <-1.3mm,0mm> at 0 1
   \put {$\scriptstyle t_4$} [rb] <-1.3mm,0mm> at 1 1
   \put {$\scriptstyle b_1$} [lb] <0.8mm,0.8mm> at 0 2
   \put {$\scriptstyle b_2$} [lb] <0.8mm,0.8mm> at 2 0
   \put {$\scriptstyle b_3$} [lb] <0.8mm,0.8mm> at 1 2
   \put {$\scriptstyle b_4$} [lb] <0.8mm,0.8mm> at 2 1

   \put {$\circ$} at 0 2
   \put {$\circ$} at 2 0
   \put {$\circ$} at 2 1
   \put {$\circ$} at 1 2
\endpicture
}} 
\caption{The $(2,2)$-box order ideal and its border}\label{fig:2-2-box}
\end{figure}

Here the bullet points correspond to the terms in the order ideal~$\OO$
and the circles to its border $\partial\OO = \{b_1,b_2,b_3,b_4\} = \{ y^2,\, x^2,\, xy^2,\, x^2y\}$.

The rim of~$\OO$ is given by the set of terms $\OOrim = \{t_2,t_3,t_4\} = 
\{ y,\, x,\, xy\}$
just inside the border, and the interior of~$\OO$ is simply $\OOint = \{t_1\} = \{ 1\}$.
Thus the tuple of rim indeterminates is 
$$
\Crim \;=\;  (c_{21}, c_{22}, c_{23}, c_{24}, c_{31}, c_{32}, c_{33}, c_{34}, c_{41}, c_{42}, 
c_{43}, c_{44})
$$
and the tuple of interior indeterminates is $\Cint = (c_{11}, c_{12}, c_{13}, c_{14})$.
\end{example}

Next we provide some additional terminology which is related to the structure of the
natural set of generators of~$I(\BO)$. It was first introduced in~\cite[Sect.~4.1]{Hui1}.

\begin{definition}\label{def-exposedIndets}
Let $\OO = \{t_1, \dots, t_\mu\}$ be an order ideal in~$\mathbb{T}^n$, and let
$\partial\OO = \{b_1, \dots, b_\nu\}$ be its border.
\begin{enumerate}
\item[(a)] For $\ell \in \{1,\dots,n\}$ and $j\in \{1,\dots,\nu\}$,
the term~$b_j$ is called {\bf $x_\ell$-exposed} if it is of 
the form $b_j = x_\ell t_i$ with $t_i\in\OO$.

\item[(b)] For $\ell\in \{1,\dots,n\}$, for $i\in \{1,\dots, \mu\}$, and for 
$j\in \{1,\dots, \nu\}$, we say that the 
indeterminate $c_{ij} \in C$ is {\bf $x_\ell$-exposed} if
$x_\ell t_i \in \partial \OO$ and if
there exists an index $j'\in \{1,\dots, \nu\}$ such that $b_j, b_{j'}$
is a next-door neighbor pair with $b_{j'} = x_\ell b_j$ or an across-the-street
neighbor pair with $x_k b_{j'} = x_\ell b_j$ for some $k\in \{1,\dots, n\}$.

\item[(c)] The subtuple of~$C$ consisting of all indeterminates $c_{ij}$
which are $x_\ell$-exposed for some $\ell\in\{1,\dots,n\}$,
is called the tuple of all {\bf exposed indeterminates} in~$C$ 
and denoted by~$\Cexp$. The indeterminates in $C \setminus \Cexp$ 
are called {\bf non-exposed}.

\end{enumerate}
\end{definition}

\begin{remark}\label{rem:intnonexp}
Notice that all interior indeterminates of~$C$ are non-exposed, i.e., that
all exposed indeterminates are rim indeterminates.
\end{remark}

The following simple example illustrates this definition.

\begin{example}{\bf (The (2,1)-Box})\label{ex-12exposed}\\
In $P=\mathbb{Q}[x,y]$, consider the order ideal $\OO = \{t_1,t_2\}$, where 
$t_1=1$ and $t_2=x$. Its border is $\partial\OO = \{b_1, b_2, b_3\}$, 
where $b_1=y$, $b_2=xy$, and $b_3 = x^2$.
\begin{figure}[ht]
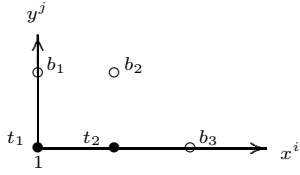

\centering{\makebox{
\beginpicture
		\setcoordinatesystem units <1cm,1cm>
		\setplotarea x from 0 to 3, y from 0 to 1.5
		\axis left /
		\axis bottom /
		\arrow <2mm> [.2,.67] from  2.5 0  to 3 0
		\arrow <2mm> [.2,.67] from  0 1  to 0 1.5
		\put {$\scriptstyle x^i$} [lt] <0.5mm,0.8mm> at 3.1 0
		\put {$\scriptstyle y^j$} [rb] <1.7mm,0.7mm> at 0 1.6

		\put {$\bullet$} at 0 0
		\put {$\bullet$} at 1 0

		\put {$\scriptstyle 1$} [lt] <-1mm,-1mm> at 0 0
		\put {$\scriptstyle t_1$} [rb] <-1.3mm,0.4mm> at 0 0
		\put {$\scriptstyle t_2$} [rb] <-1.3mm,0.4mm> at 1 0
		\put {$\scriptstyle b_3$} [lb] <0.8mm,0.4mm> at 2 0
		\put {$\scriptstyle b_2$} [rb] <4.3mm,0mm> at 1 1
		\put {$\scriptstyle b_1$} [lb] <0.8mm,0mm> at 0 1

  		\put {$\circ$} at 0 1
		\put {$\circ$} at 1 1
		\put {$\circ$} at 2 0
\endpicture
}} 
\caption{The $(2,1)$-box order ideal and its border}\label{fig:2-1-box}%
\end{figure}

There are two neighbor pairs, namely the next-door pair $b_2 = x b_1$ 
and the across-the-rim pair $y b_3 = x b_2$. As they involve 
$x b_1$ and $x b_2$, we get $j\in \{1,2\}$.
Consequently, the $x$-exposed indeterminates $c_{ij}$ in~$C$
satisfy $i\in \{ 2\}$  and $j\in \{1,2\}$, i.e., they are $c_{21}$ and~$c_{22}$.
Similarly, we see that the $y$-exposed indeterminates satisfy $i\in \{1,2 \}$ 
and $j \in \{3\}$, i.e., they are $c_{13}$ and~$c_{23}$.

Altogether, we obtain $\Cexp = (c_{13}, c_{21}, c_{22}, c_{23})$.
\end{example}

The next proposition provides two characterizations of exposed
indeterminates which motivate the above definition.

\begin{proposition}{\bf (Characterizing Exposed Indeterminates)}\label{prop:CharExposed}\\
Let $\OO = \{t_1,\dots,t_\mu\}$ be an order ideal in~$\mathbb{T}^n$, let $\partial\OO = 
\{b_1,\dots,b_\nu\}$ be its border, and
let $G = \{g_1,\dots,  g_\nu\}$ be the generic $\OO$-border 
prebasis, where $g_j \;= \; b_j - \sum_{i=1}^\mu c_{ij} t_i$ for $j=1\dots, \nu$.
Then, for all $i\in \{1,\dots, \mu\}$ and $j\in \{1,\dots, \nu\}$, 
the following conditions are equivalent.
\begin{enumerate}
\item[(a)] There exists an index $\ell \in \{1,\dots, n\}$ such that
the indeterminate $c_{ij} \in C$ is $x_\ell$-exposed.

\item[(b)] There exists an index $j'\in \{1,\dots, \nu\}$ such that $b_j, b_{j'}$
is a next-door neighbor pair with $b_{j'} = x_\ell b_{j}$ or an across-the-street
neighbor pair with $x_k b_{j'} = x_\ell b_j$ for some $k\in \{1,\dots, n\}$, and
such that if we write 
$$
x_\ell g_j \;=\; x_\ell b_j  + c_{i_1 j} g_{j_1} + \cdots + c_{i_r j} g_{j_r} 
+ \tsum_{\lambda=1}^\mu f_\lambda t_\lambda
$$
where $i_1,\dots,  i_r\in \{1,\dots, \mu\}$, where $j_1,\dots, j_r\in\{1,\dots, \nu\}$, 
and where $f_\lambda \in K[C]$, then we have $i \in \{i_1, \dots,  i_r \}$.

\item[(c)] There exists an index $j'\in \{1,\dots, \nu\}$ such that $b_j, b_{j'}$
is a next-door neighbor pair with $b_{j'} = x_\ell b_j$ or an across-the-street
neighbor pair with $x_k b_{j'} = x_\ell b_j$ and $k\in \{1,\dots, n\}$, and
such that the lifting of the neighbor syzygy $e_{j'} - x_\ell e_j$ resp.\ 
$x_k e_{j'} - x_\ell e_j$ is of the form 
$$
x_k^\epsilon e_{j'} - x_\ell e_j - \tsum_{\kappa=1}^\nu \bar{f}_\kappa  e_{\kappa}
\in \Syz_{B_\OO}( \bar{g}_1, \dots,  \bar{g}_\nu )
$$
where $\epsilon\in \{0,1\}$ and $\bar{f}_\kappa \in B_\OO$ are the residue classes 
of monomials or binomials in $K[C]$ such that one of their supports contains $c_{ij}$.

\end{enumerate}
\end{proposition}

\begin{proof}
To show (a)$\Rightarrow$(b), we write $x_\ell t_i = b_{j''}$ with 
$j'' \in \{1,\dots, \nu\}$. Clearly, there exists a neighbor pair $b_j, b_{j'}$
such that $x_k^\epsilon b_{j'} = x_\ell b_j$ for some $j, j' \in \{1,\dots, \nu\}$
and $k\in \{1,\dots, n\}$ and $\epsilon \in \{0,1\}$. (For instance, we can take the
pure $x_\ell$-power $b_j = x_\ell^\kappa$ in~$\partial\OO$.) Then we have
$$
x_\ell g_j \;=\; x_\ell b_j - \tsum_{\lambda=1}^\mu  c_{\lambda j} x_\ell t_\lambda
$$
and if we denote those terms $x_\ell t_{\lambda}$ which are in $\partial\OO$ by
$x_\ell t_{i_1} = b_{j_1}$, $\dots$, $x_\ell t_{i_r} = b_{j_r}$ then
the fact that $x_\ell t_i \in \partial\OO$ shows $i \in \{i_1,\dots,  i_r\}$.

Conversely, the condition $f_\lambda \in K[C]$ implies that the terms in the
support of $x_\ell (g_j-b_j)$ 
which are contained in $K[C] \cdot \partial \OO$ are of the form $c_{i_\kappa j} b_{j_\kappa}$
with $\kappa \in \{1,\dots, r\}$. Therefore $b_{j_\kappa} = x_\ell t_{i_\kappa}$
for $\kappa = 1,\dots, r$ and $c_{ij} \in \{ c_{i_1 j_1}, \dots,  c_{i_r j_r} \}$ 
is $x_\ell$-exposed.

The equivalence of~(b) and~(c) follows from the fact that the representation of
$x_k^\epsilon g_{j'} - x_\ell g_j$ which follows from~(b) is exactly the lifting
of the neighbor syzygy $x_k^\epsilon e_{j'} - x_\ell e_j$ in 
$\Syz_{B_\OO}( \bar{g}_1, \dots,  \bar{g}_\nu )$ by~\cite{KR3}, Prop. 4.1.
\end{proof}

\bigskip\bigbreak
%
%

\section{MaxDeg Border Basis Schemes}
\label{sec3}

Our main goal in this paper is to find good re-embeddings of certain border basis schemes.
For this purpose, we intend to use the technique of $Z$-separating re-embeddings which
has been developed in~\cite{KLR1,KLR2,KLR3,KR5}. Let us recall the basic construction
using the current setting.

Let $K$ be a field, let $P=K[x_1,\dots,x_n]$, let $\OO=\{t_1,\dots,t_\mu\}$
be an order ideal in~$\mathbb{T}^n$, let $\partial\OO = \{b_1,\dots,b_\nu\}$
be the border of~$\OO$, let $C=(c_{11},\dots,c_{\mu\nu})$ be a tuple of
indeterminates, and let $I(\BO) \subseteq K[C]$ be the defining ideal of the $\OO$-border basis scheme.

\begin{definition}{\bf (Separating and Coherently Separating Tuples)}\label{def-sepindets}\,\\
Let $Z=(z_1,\dots, z_s)$ be a tuple of distinct indeterminates in $C=(c_{11},\dots,c_{\mu\nu})$.
\begin{enumerate}
\item[(a)] Let $i\in \{1, \dots, s\}$. A polynomial $f\in K[C]$ is called {\bf $z_i$-separating}
if it is of the form $f = z_i - h$ with a polynomial~$h$ such that $z_i$ divides
no term in $\Supp(h)$.

\item[(b)] A tuple of polynomials $F=(f_1,\dots, f_s) \in K[C]^s$ is called
{\bf $Z$-separating} if there exists a term ordering~$\sigma$ such that
$\LT_\sigma(f_i)=z_i$ for $i=1, \dots, s$.

\item[(c)] A tuple of polynomials $F=(f_1,\dots, f_s) \in K[C]^s$ is called
{\bf coherently $Z$-separating} if the polynomials~$f_i$ are of the form $f_i = z_i - h_i$
with polynomials~$h_i$ such that $z_i$ divides no term in $\Supp(h_j)$
for $i,j \in \{1,\dots, s\}$.

\item[(d)] Let $I$ be an ideal in~$K[C]$. The tuple~$Z$ is called a  
{\bf separating tuple of indeterminates} for~$I$ if~$I$ contains a $Z$-separating 
tuple of polynomials.

\end{enumerate}
\end{definition}

The usefulness of coherently $Z$-separating tuples of polynomials derives from the fact that
they allow us to perform re-embeddings of ideals in~$K[C]$ as follows.

\begin{proposition}{\bf ($Z$-Separating Re-Embeddings)}\label{prop-ZsepReemb}\\
Let $I\subseteq K[C]$ be an ideal which is 
generated by a tuple of polynomials $G=\langle g_1,\dots, g_r)$, let $Z=(z_1,\dots, z_s)$
be a tuple of distinct indeterminates in~$C$, and let $F=(f_1,\dots, f_s)$ be a coherently
$Z$-separating tuple of polynomials in~$I$. For $i=1, \dots, s$, write
$f_i = z_i - h_i$ with $h_i\in K[C]$.
For $i=1,\dots, s$ and $j=1, \dots, r$, replace each occurrence of~$z_i$ in~$g_j$ by~$h_i$. 
The resulting polynomial $\hat{g}_j$ is said to be obtained by {\bf rewriting~$g_j$ using~$F$}.
\begin{enumerate}
\item[(a)] The tuple $\widehat{G} = (\hat{g}_1, \dots, \hat{g}_r)$ obtained by rewriting~$G$ 
using~$F$ is a system of generators of the elimination ideal $I\cap K[C\setminus Z]$.

\item[(b)] The $K$-algebra homomorphism
$$
\Phi:\; K[C]/I \longrightarrow K[C\setminus Z] / (I\cap K[C\setminus Z]) 
$$
defined by $z_i \mapsto h_i$ for $i=1,\dots,  s$ and $c_{kj} \mapsto c_{kj}$ for
$c_{kj} \notin Z$ is an isomorphism. It is called a {\bf $Z$-separating re-embedding} of~$I$
(or of $K[C]/I$).

\end{enumerate}
\end{proposition}

A proof of this proposition is contained in~\cite[Thm.~2.10]{KLR1}. 
Notice that performing the isomorphism~$\Phi$ amounts to substituting the polynomials~$h_i$
for the indeterminates~$z_i$ finitely often.

The ideal $I \cap K[C\setminus Z]$ defines an affine scheme which is isomorphic
to the original scheme $\Spec(K[C]/I)$ but embedded in a lower dimensional affine space. 
This makes it easier to study this scheme, and sometimes we may even be able to prove that
it is actually isomorphic to an affine space (see, for instance, Propositions~\ref{prop-L-Shape-Zsep} 
and~\ref{prop-theLshape}).

\medskip

In~\cite{KR5} we showed that the method of $Z$-separating re-embeddings
is particularly powerful when the ideal~$I$ is homogeneous with respect to a non-negative grading. 
Since the ideals $I(\BO)$ are homogeneous with respect to the total arrow grading,
we are interested in knowing when this grading is non-negative. 
This property can be characterized as follows (see also~\cite[Sect.~2]{KR3}).

\begin{proposition}\label{prop-CharMaxDeg}
Let $\OO = \{t_1,\dots,  t_\mu\}$ be an order ideal with border $\partial\OO = 
\{b_1,\dots,  b_\nu\}$, and let $I(\BO) \subseteq K[C]$ be the ideal defining the $\OO$-border
basis scheme $\BO$. Then the following conditions are equivalent.
\begin{enumerate}
\item[(a)] The order ideal $\OO$ has a {\bf MaxDeg border}, i.e., we have
$\deg(b_j) \ge \deg(t_i)$ for $i=1, \dots, \mu$ and $j=1, \dots, \nu$.

\item[(b)] The total arrow grading on~$K[C]$ is non-negative.

\item[(c)] The order ideal $\OO$ has a {\bf generic Hilbert function}, i.e., 
there exists a degree $d\ge 0$ such that $\# \OO_i = \HF_P(i)$ for $i<d$,
such that $\#\OO_d >0$, and such that $\# \OO_i = 0$ for $i>d$.

\end{enumerate}
\end{proposition}

\begin{proof}
The equivalence of~(a) and~(b) follows immediately from the definition.
The equivalence of~(a) and~(c) was shown in~\cite[Prop.~5.8]{KLR0}.
\end{proof}

In order to apply the results of~\cite{KR5}, we need one further
ingredient. Namely, in~\cite{KR5}, it is assumed that the ideal~$I$
which we want to re-embed satisfies $I \cap P_0 = \{0\}$, i.e.\ it does
not contain any non-zero homogeneous elements of degree zero.
This is the content of the following proposition which is based on
a result in~\cite{KR3}.

\begin{proposition}\label{prop-MaxDegBasis}
Let $\OO$ be an order ideal in~$\mathbb{T}^n$ which has a MaxDeg border.
Recall that $\Cnull$ is the subtuple of~$C$ consisting of all indeterminates
of total arrow degree zero. Then the following claims hold.
\begin{enumerate}
\item[(a)] We have $I(\BO) \cap K[\Cnull] = \{0\}$.

\item[(b)] The ring $K[\Cnull]$ is the homogeneous component of degree
zero of the non-negatively graded $K$-algebra $B_\OO = K[C]/I(\BO)$.
\end{enumerate}
\end{proposition}

\begin{proof}
By assumption, we have $\deg_W(c_{ij}) \ge 0$ for all $i,j$, where~$W$
defines the total arrow grading. For $i=1,\dots,n$, let $\overline{\mathcal{A}}_i$
be the matrix obtained from the generic multiplication matrix $\mathcal{A}_i$ by 
setting $c_{ij} \mapsto 0$ whenever $\deg_W(c_{ij})>0$.
The commutators of the matrices $\overline{\mathcal{A}}_i$ define the homogeneous
border basis scheme $\BOhom$ (see~\cite[Sect.~5]{KR3}). It is shown 
in~\cite[Thm.~5.3]{KR3} that the matrices $\overline{\mathcal{A}}_i$ commute.
Therefore no commutator $\mathcal{A}_i \mathcal{A}_j - \mathcal{A}_j \mathcal{A}_i$
contains a non-zero homogeneous polynomial of total arrow degree zero. Consequently, all
generators of~$I(\BO)$ have positive arrow degrees and claim~(a) follows.

Claim~(b) is an immediate consequence of~(a).
\end{proof}

Notice that it is easy to count the number of indeterminates $\# \Cnull$
in the polynomial ring $K[\Cnull]$. Since we have $\deg(b_j) \ge \deg(t_i)$
for all $i,j$, and since~$\OO$ has a generic Hilbert function,
the only way to get $\deg_W(c_{ij})= 0$ is from
$\deg(t_i)=\deg(b_j)=d$, where $d=\max\{ \deg(t_i) \mid i=1,\dots,\mu\}$.
Therefore we have 
$$
\# \Cnull \;=\; \# \{ i\in \{1,\dots,\mu\} \mid \deg(t_i)=d\} \;\cdot\;
\# \{ j\in \{1,\dots,\nu\} \mid \deg(b_j) = d\}
$$

For a MaxDeg order ideal $\OO$, the ring $(B_\OO)_0$ is the coordinate ring of 
the {\bf homogeneous border basis scheme} $\BOhom$. 
This is the closed subscheme of~$\BO$ which parametrizes all
$\OO$-border bases which are homogeneous with respect to the standard grading.
The above proposition can be interpreted by saying that $\BOhom$ is an affine cell, 
more precisely that there is an isomorphism of schemes 
$\BOhom \cong \mathbb{A}^{\gamma}$, where $\gamma = \# \Cnull$.
Moreover, the inclusion map $(B_\OO)_0 \longrightarrow B_\OO$ corresponds to a
dominant morphism of affine schemes $\Theta:\; \BO \longrightarrow \mathbb{A}^\gamma$.
In fact, this morphism is surjective, because the coordinate rings of the 
fibers are positively graded.

The next proposition collects some properties of MaxDeg border basis schemes
which follow from the results of~\cite{KR5}.
Recall that the {\bf linear part} $\lin(g)$ of a polynomial $g\in K[C]$ is its homogeneous component of
(standard) degree one, and that the linear parts of the natural generators of~$I(\BO)$ can be described
explicitly (see~\cite[Cor.~2.8]{KSL} and~\cite[Prop.~6.3]{KLR3}).

\begin{proposition}{\bf (Re-Embeddings of MaxDeg Border Basis Schemes)}\label{prop-MaxDegReEmbed}\\
Let $\OO$ be an order ideal in~$\mathbb{T}^n$ which has a MaxDeg border, and let
$\{g_1, \dots, g_r\}$ be a set of generators of~$I(\BO)$ consisting of polynomials which are
homogeneous with respect to the total arrow grading.
\begin{enumerate}
\item[(a)] Let $Z$ be a tuple of indeterminates which is separating for~$I(\BO)$,
i.e., for which a $Z$-separating re-embedding of~$I(\BO)$ exists. Then we have $Z \subseteq \Cplus$,
i.e., $Z$ consists of indeterminates of positive total arrow degree.

\item[(b)] A tuple of indeterminates $Z=(z_1,\dots,z_s)$ is separating for~$I(\BO)$
if and only if $z_1,\dots,z_s \in \langle \lin(g_1),\dots,\lin(g_r) \rangle_{K[\Cnull]}$.

\item[(c)] Suppose that $K$ is perfect and $\BO$ is a regular scheme. Then there exists a $K$-algebra
isomorphism $B_\OO \cong K[\Cnull][\widehat{X}]$ with a tuple of indeterminates $\widehat{X}$.
In other words, the scheme $\BO$ is an {\bf affine cell}.

\end{enumerate}
\end{proposition}

\begin{proof}
Claim (a) follows from~\cite[Rem.~3.1.b]{KR5}, and~(b) is a consequence of~\cite[Thm.~3.3.c]{KR5}.
The isomorphism in~(c) requires us to combine a $Z$-separating re-embedding of~$I(\BO)$
with an isomorphism based on Unimodular Matrix Problem (see~\cite[Thm.~6.7]{KR5}).
\end{proof}

Notice that the isomorphism in part~(c) of this proposition can be computed 
effectively (see~\cite[Alg.~6.8]{KR5}). In the next section we carry out one such 
computation explicitly. An important class of examples is given as follows.

\begin{corollary}{\bf (Planar MaxDeg Border Basis Schemes)}\label{cor-bivariateaffinecells}\\
Let $K$ be a perfect field. Then every planar MaxDeg  border basis scheme 
defined over $K$ is an affine cell. More precisely, there exists a homogeneous
$K$-algebra isomorphism $B_\OO \cong K[\widetilde{X}]$, where $\widetilde{X}$ 
is a tuple of $\mu\nu$ indeterminates.
\end{corollary}

\begin{proof}
By~\cite[Thm.~2.4]{Fog}, the Hilbert scheme $\Hilb^\mu(\mathbb{A}^2)$ is regular.
Since the border basis scheme $\BO$ is embedded as an open subscheme of
$\Hilb^\mu(\mathbb{A}^2)$ by~\cite[Lemma 7]{Hui3}, it is also regular,
and part~(c) of the proposition applies.
\end{proof}

The inclusion $K[\Cnull] = (B_\OO)_0 
\hookrightarrow B_\OO$ allows us to view $\BO$ as a family of schemes, parametrized by
the affine space $\Spec(K[\Cnull]) \cong \mathbb{A}^\gamma$. This point of view 
makes it clear that we have a special case of~\cite[Sect.~5.4]{KR5} and can apply
the results shown there about the fibers of the family $\theta:\, \BO \longrightarrow
\mathbb{A}^\gamma$. For instance, we obtain the following insights.

\begin{proposition}{\bf (The MaxDeg Border Basis Family)}\label{prop-MaxDegFamily}\\
Let $\OO$ be an order ideal in~$\mathbb{T}^n$ which has a MaxDeg border,
let $\gamma = \# \Cnull$, let $\Gamma \in K^\gamma$, and let
$\m_\Gamma \subseteq K[\Cnull]$ be the maximal ideal corresponding to~$\Gamma$.
\begin{enumerate}
\item[(a)] The ideal of the fiber of~$\theta$ over~$\Gamma$ is
$I(\BO)_\Gamma = I(\BO) \cdot (K[\Cnull]/\m_\Gamma)[\Cplus] \subseteq K[\Cplus]$
and there exists a tuple of indeterminates~$Z$ in~$\Cplus$ which yields a 
separating re-embedding
$$
\Phi_\Gamma:\; K[\Cplus] / I(\BO)_\Gamma \;\longrightarrow\; 
K[\Cplus \setminus Z] / I(\BO)_\Gamma \cap K[\Cplus \setminus Z]
$$
This re-embedding is {\bf optimal} in the sense that it reaches 
the minimal embedding dimension.

\item[(b)] If the localization of $K[\Cplus]/I(\BO)_\Gamma$ at its irrelevant maximal ideal
is regular, then the map $\Phi_\Gamma$ defines an isomorphism between the coordinate ring of the fiber
$K[\Cplus]/I(\BO)_\Gamma$ and a polynomial ring. In other words, in this case the fiber
of~$\theta$ over~$\Gamma$ is an affine space.

\end{enumerate}
\end{proposition}

\begin{proof} Claim (a) is a consequence of~\cite[Prop.~5.4.b]{KR5}, and
claim (b) is shown in~\cite[Thm.~5.9.d]{KR5}.
\end{proof}

Notice that~\cite{KR5} also provides similar results for the generic fiber of~$\phi$.
If we drop the assumption that~$\OO$ has a MaxDeg border and set
$c_{ij}\mapsto 0$ for all indeterminates in~$C$ with negative total arrow degree
$\deg_W(c_{ij}) <0$, we get a closed subscheme of~$\BO$ which
is non-negatively graded, as well. This subscheme $\BOdf$ is called
the {\bf degree filtered border basis scheme} and was studied in~\cite[Sect.~5]{KLR0}.
Its coordinate ring is non-negatively graded, and it has the coordinate ring of
$\BOhom$ as its homogeneous component of degree zero.

However, unfortunately, as the following examples show, the defining ideal
$I(\BOdf)$ does, in general, intersect the ring $K[\Cnull]$ non-trivially. Hence
the homogeneous border basis scheme is not equal to $\Spec(K[\Cnull])$ and the
results of~\cite{KR5} cannot be applied directly to this situation.

\begin{example}\label{ex-motaffine1}
Consider the polynomial ring $P=\QQ[x,y]$ and the order ideal $\OO = \{1, y, x, y^2, xy, y^3 \}$. 
Here we have $\partial\OO = \{x^2,  xy^2,  x^2y,  y^4,  xy^3\}$ (see Figure~\ref{fig:NonMaxDeg-1}).

\begin{figure}[ht]
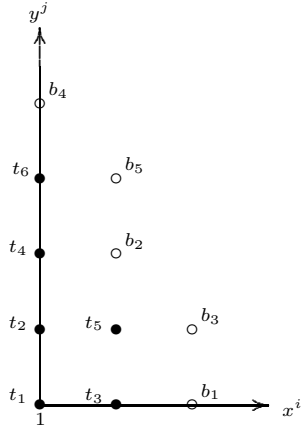

\centering{\makebox{
\beginpicture
		\setcoordinatesystem units <1cm,1cm>
		\setplotarea x from 0 to 3, y from 0 to 4.5
		\axis left /
		\axis bottom /
		\arrow <2mm> [.2,.67] from  2.5 0  to 3 0
		\arrow <2mm> [.2,.67] from  0 4.5  to 0 5
		\put {$\scriptstyle x^i$} [lt] <0.5mm,0.8mm> at 3.1 0
		\put {$\scriptstyle y^j$} [rb] <1.7mm,0.7mm> at 0 5

		\put {$\bullet$} at 0 0
		\put {$\bullet$} at 1 0
		\put {$\bullet$} at 0 1
		\put {$\bullet$} at 1 1
		\put {$\bullet$} at 0 2
		\put {$\bullet$} at 0 3

		\put {$\scriptstyle 1$} [lt] <-1mm,-1mm> at 0 0
		\put {$\scriptstyle t_1$} [rb] <-1.3mm,0.4mm> at 0 0
		\put {$\scriptstyle t_3$} [rb] <-1.3mm,0.4mm> at 1 0
		\put {$\scriptstyle t_2$} [rb] <-1.3mm,0mm> at 0 1
		\put {$\scriptstyle t_4$} [rb] <-1.3mm,0mm> at 0 2
		\put {$\scriptstyle t_5$} [rb] <-1.3mm,0mm> at 1 1
		\put {$\scriptstyle t_6$} [lb] <-4mm,0.4mm> at 0 3
		
		\put {$\scriptstyle b_1$} [lb] <0.8mm,0.8mm> at 2 0
		\put {$\scriptstyle b_2$} [lb] <0.8mm,0.8mm> at 1 2
		\put {$\scriptstyle b_3$} [lb] <0.8mm,0.8mm> at 2 1
	    \put {$\scriptstyle b_4$} [lb] <0.8mm,0.8mm> at 0 4
	    \put {$\scriptstyle b_5$} [lb] <0.8mm,0.8mm> at 1 3
	
		\put {$\circ$} at 2 0
		\put {$\circ$} at 2 1
		\put {$\circ$} at 1 2
		\put {$\circ$} at 1 3
		\put {$\circ$} at 0 4
\endpicture
}} 
\caption{An order ideal without MaxDeg border}\label{fig:NonMaxDeg-1}
\end{figure}

Since $\deg(y^3) > \deg(x^2)$,
the order ideal~$\OO$ does not have a MaxDeg border. In order to define the degree filtered
border basis scheme $\BOdf$, we have to set $c_{61} \mapsto 0$ in $I(\BO)$.
In other words, we have $I(\BOdf) = I(\BO) + \langle c_{61}\rangle \subseteq K[C]$.

In this case it is easy to check that $I(\BOdf) \cap K[\Cnull] = \langle
c_{41} -c_{63}  +c_{51}c_{62}\rangle$, so the hypotheses of~\cite{KR5} are not satisfied.
Notice that the homogeneous border basis scheme $\BOhom$ has the
coordinate ring $\QQ[c_{41}, c_{51}, c_{62}, c_{63}] / \langle c_{41} -c_{63}  +c_{51}c_{62}\rangle$
here. It is isomorphic to $\mathbb{A}^3$, and using the substitution $c_{63} \mapsto
c_{41} + c_{51} c_{63}$ we could still transform everything to the setting of~\cite{KR5}.
\end{example}

\begin{example}\label{ex-notaffine2}
Consider the polynomial ring $P=\QQ[x,y]$ and the order ideal $\OO=\{1, y, x, y^2, x^2, y^3 \}$. 
Here we have $\partial\OO = \{xy,  xy^2,  x^2y,  x^3,  y^4,  xy^3\}$ (see Figure~\ref{fig:NonMaxDeg-2}). 

\begin{figure}[ht]
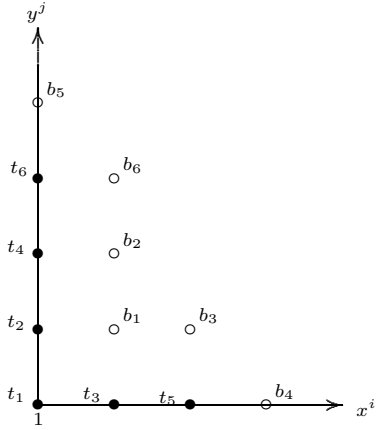

\centering{\makebox{
\beginpicture
		\setcoordinatesystem units <1cm,1cm>
		\setplotarea x from 0 to 4, y from 0 to 4.5
		\axis left /
		\axis bottom /
		\arrow <2mm> [.2,.67] from  3.5 0  to 4 0
		\arrow <2mm> [.2,.67] from  0 4.5  to 0 5
		\put {$\scriptstyle x^i$} [lt] <0.5mm,0.8mm> at 4.1 0
		\put {$\scriptstyle y^j$} [rb] <1.7mm,0.7mm> at 0 5

		\put {$\bullet$} at 0 0
		\put {$\bullet$} at 1 0
		\put {$\bullet$} at 0 1
		\put {$\bullet$} at 2 0
		\put {$\bullet$} at 0 2
		\put {$\bullet$} at 0 3

		\put {$\scriptstyle 1$} [lt] <-1mm,-1mm> at 0 0
		\put {$\scriptstyle t_1$} [rb] <-1.3mm,0.4mm> at 0 0
		\put {$\scriptstyle t_3$} [rb] <-1.3mm,0.4mm> at 1 0
		\put {$\scriptstyle t_2$} [rb] <-1.3mm,0mm> at 0 1
		\put {$\scriptstyle t_4$} [rb] <-1.3mm,0mm> at 0 2
		\put {$\scriptstyle t_5$} [rb] <-1.3mm,0mm> at 2 0
		\put {$\scriptstyle t_6$} [lb] <-4mm,0.4mm> at 0 3
		
		\put {$\scriptstyle b_1$} [lb] <0.8mm,0.8mm> at 1 1
		\put {$\scriptstyle b_2$} [lb] <0.8mm,0.8mm> at 1 2
		\put {$\scriptstyle b_3$} [lb] <0.8mm,0.8mm> at 2 1
	    \put {$\scriptstyle b_4$} [lb] <0.8mm,0.8mm> at 3 0
	    \put {$\scriptstyle b_5$} [lb] <0.8mm,0.8mm> at 0 4
	    \put {$\scriptstyle b_6$} [lb] <0.8mm,0.8mm> at 1 3
	
		\put {$\circ$} at 1 1
		\put {$\circ$} at 2 1
		\put {$\circ$} at 1 2
		\put {$\circ$} at 1 3
		\put {$\circ$} at 0 4
		\put {$\circ$} at 3 0
\endpicture
}} 
\caption{Another order ideal without MaxDeg border}\label{fig:NonMaxDeg-2}
\end{figure}

Since $\deg(y^3) > \deg(xy)$,
the order ideal~$\OO$ does not have a MaxDeg border. The degree filtered border basis scheme is defined
by the ideal $I(\BOdf) = I(\BO) + \langle c_{61} \rangle$.

In this case we calculate $I(\BOdf) \cap K[\Cnull] = \langle f_1,\, f_2\rangle$,
where $f_1=c_{41} -c_{62} +c_{51}c_{63}$ and $f_2=c_{63} -c_{41}c_{62} - c_{51}c_{64}$.
Again this means that the results of~\cite{KR5} cannot be applied directly to~$\BOdf$.

The coordinate ring of $\BOhom$ is $\QQ[c_{41}, c_{51}, c_{62}, c_{63}, c_{64}]/\langle f_1,f_2\rangle$.
As there is no pair of indeterminates~$Z$ which is separating for the ideal $\langle f_1,f_2\rangle$,
it is not clear whether $\BOhom$ is isomorphic to~$\mathbb{A}^3$.
\end{example}

\bigskip\bigbreak
%
%

\section{The L-Shape Border Basis Scheme}
\label{sec4}

In this section we apply the techniques developed above to a case which, in spite
of its apparent simplicity, proved to be quite enigmatic for some time.
We work in the bivariate polynomial ring $P=K[x,y]$ over a computable perfect field~$K$
and consider the following border basis scheme.

\begin{definition}\label{def:L-shape}
The order ideal $\OO = \{1,\, y,\, x,\, y^2,\, x^2\}$ in~$\mathbb{T}^2$
is called the {\bf L-shape order ideal}.
\end{definition}

In this section, we always let $\OO=\{t_1,t_2,t_3,t_4,t_5\}$ be the L-shape order ideal, 
where $t_1=1$, $t_2=y$, $t_3=x$, $t_4=y^2$, and
$t_5=x^2$. The border of~$\OO$ is given by $\partial\OO = \{b_1,\dots,b_5\}$,
where $b_1=xy$, $b_2=y^3$, $b_3=xy^2$, $b_4=x^2y$, and $b_5=x^3$.
The following figure depicts the L-shape order ideal and its border.
\begin{figure}[ht]
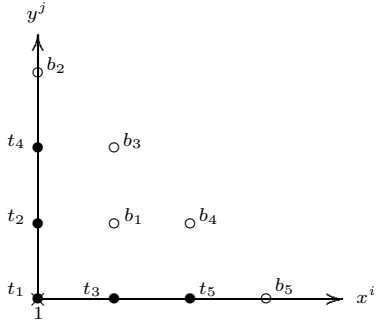

\centering{\makebox{
\beginpicture
		\setcoordinatesystem units <1cm,1cm>
		\setplotarea x from 0 to 4, y from 0 to 3.5
		\axis left /
		\axis bottom /
		\arrow <2mm> [.2,.67] from  3.5 0  to 4 0
		\arrow <2mm> [.2,.67] from  0 3  to 0 3.5
		\put {$\scriptstyle x^i$} [lt] <0.5mm,0.8mm> at 4.1 0.1
		\put {$\scriptstyle y^j$} [rb] <1.7mm,0.7mm> at 0 3.6
		\put {$\bullet$} at 0 0
		\put {$\bullet$} at 1 0
		\put {$\bullet$} at 0 1
		\put {$\bullet$} at 0 2
		\put {$\bullet$} at 2 0
		\put {$\scriptstyle 1$} [lt] <-1mm,-1mm> at 0 0
		\put {$\scriptstyle t_1$} [rb] <-1.3mm,0.4mm> at 0 0
		\put {$\scriptstyle t_3$} [rb] <-1.3mm,0.4mm> at 1 0
		\put {$\scriptstyle t_2$} [rb] <-1.3mm,0mm> at 0 1
		\put {$\scriptstyle t_4$} [rb] <-1.3mm,0mm> at 0 2
		\put {$\scriptstyle t_5$} [lb] <0.8mm,0.4mm> at 2 0
		\put {$\scriptstyle b_1$} [rb] <4.3mm,0mm> at 1 1
		\put {$\scriptstyle b_2$} [lb] <0.8mm,0mm> at 0 3
		\put {$\scriptstyle b_3$} [lb] <0.8mm,0mm> at 1 2
		\put {$\scriptstyle b_4$} [lb] <0.8mm,0mm> at 2 1
		\put {$\scriptstyle b_5$} [lb] <0.8mm,0.7mm> at 3 0
		\put {$\times$} at 0 0
		\put {$\circ$} at 0 3
		\put {$\circ$} at 2 1
		\put {$\circ$} at 1 2
		\put {$\circ$} at 1 1
		\put {$\circ$} at 3 0
\endpicture
}} 
\caption{The L-shape order ideal and its border}\label{fig:L-shape}
\end{figure}

\begin{definition}
For the L-shape order ideal~$\OO$, the corresponding border basis scheme~$\BO$
is called the {\bf L-shape border basis scheme}.
\end{definition}

Since $\mu=\nu=5$, the L-shape border basis scheme $\BO$ is naturally embedded into
the affine space $\mathbb{A}^{25}$. Its vanishing ideal is generated by the 20 nonzero entries of the
commutator $\mathcal{A}_1 \mathcal{A}_2 - \mathcal{A}_2 \mathcal{A}_1$ of the two
generic multiplication matrices. As we know, the ideal $I(\BO)$ is homogeneous
with respect to the total arrow grading which assigns to the tuple of indeterminates
$C=(c_{11}, c_{12}, \dots, c_{55})$ the weights given by
$$
W = (2,  3,  3,  3,  3,  1,  2,  2,  2,  2,  1,  2,  2,  2,  2,  0,  1,  1,  1,  1,  0,  1,  1,  1,  1)
$$
In particular, we see that $\OO$ has a MaxDeg border and that the total arrow grading is
non-negative, so that we can apply the results of the preceding section.
We begin by using the technique of $Z$-separating re-embeddings.

\begin{proposition}\label{prop-L-Shape-Zsep}
The tuple of indeterminates
$$
Z \;=\; (c_{11},\; c_{12},\; c_{13},\; c_{14},\; c_{15},\; c_{23},\; 
     c_{24},\; c_{25}, \; c_{31}, \; c_{32},\; c_{34},\; c_{44},\; c_{53})
$$
defines a {\it best} separating re-embedding of~$I(\BO)$, i.e., a separating
re-embedding for which $\#Z$ is as large as possible.
This re-embedding is a $K$-algebra isomorphism 
$$
\Phi:\; B_\OO \;\longrightarrow\; K[C\setminus Z]/ \langle f_1,\, f_2\rangle
$$
where
{\small
$$
\begin{array}{lcl}
f_1 &=&   c_{21}c_{41}^2c_{51}^2  +c_{41}^2c_{43}c_{51}^2 +c_{41}c_{45}c_{51}^3 +c_{41}^2c_{42}c_{51} 
   -c_{41}^3c_{52} -c_{41}^2c_{51}c_{54}\\
&& +c_{21}c_{41}c_{51} +c_{45}c_{51}^2 -c_{41}c_{51}c_{55} +c_{41}c_{42} +c_{41}c_{54} +c_{21} -c_{43} 
\medskip
\cr
f_2 &=&  -c_{21}^2c_{41}^3c_{51}^5 -2c_{21}c_{41}^3c_{43}c_{51}^5 -c_{41}^3c_{43}^2c_{51}^5 -2c_{21}c_{41}^2c_{45}c_{51}^6 
   -2c_{41}^2c_{43}c_{45}c_{51}^6 \\
&&-c_{41}c_{45}^2c_{51}^7 -c_{21}c_{41}^3c_{42}c_{51}^4 -c_{41}^3c_{42}c_{43}c_{51}^4 -c_{41}^2c_{42}c_{45}c_{51}^5 
   +c_{21}c_{41}^3c_{51}^4c_{54} \\
&&+c_{41}^3c_{43}c_{51}^4c_{54} +c_{41}^2c_{45}c_{51}^5c_{54} -c_{41}^4c_{42}c_{51}^2c_{52} +c_{41}^5c_{51}c_{52}^2 
   +c_{41}^4c_{51}^2c_{52}c_{54}\\
&& -c_{41}^2c_{42}c_{43}c_{51}^3 -c_{21}c_{41}^3c_{51}^2c_{52} +c_{41}^3c_{43}c_{51}^2c_{52} -c_{41}^2c_{45}c_{51}^3c_{52}\\
&& -2c_{21}c_{41}^2c_{51}^3c_{54} -c_{41}^2c_{43}c_{51}^3c_{54} -2c_{41}c_{45}c_{51}^4c_{54} -c_{41}^2c_{42}c_{51}^3c_{55}\\
&& +2c_{41}^3c_{51}^2c_{52}c_{55} +c_{41}^2c_{51}^3c_{54}c_{55} +2c_{21}c_{41}c_{43}c_{51}^3 +3c_{41}c_{43}^2c_{51}^3 
   +2c_{43}c_{45}c_{51}^4 \\
&&-c_{41}^4c_{52}^2 -2c_{41}^3c_{51}c_{52}c_{54}
+2c_{41}c_{43}c_{51}^3c_{55} +c_{41}c_{51}^3c_{55}^2 +c_{41}c_{42}c_{43}c_{51}^2\\
&& -c_{42}c_{45}c_{51}^3 +c_{21}c_{41}^2c_{51}c_{52} +3c_{41}c_{45}c_{51}^2c_{52} -3c_{41}c_{43}c_{51}^2c_{54} +c_{45}c_{51}^3c_{54}\\
&& +c_{41}c_{42}c_{51}^2c_{55} -3c_{41}^2c_{51}c_{52}c_{55} -3c_{41}c_{51}^2c_{54}c_{55} +c_{22}c_{41}c_{51} -2c_{33}c_{41}c_{51}\\
&& +c_{21}^2c_{51}^2 -c_{35}c_{51}^2 -c_{43}^2c_{51}^2 +c_{41}^2c_{42}c_{52} +c_{41}^2c_{52}c_{54} +c_{41}c_{51}c_{54}^2 \\
&&-2c_{43}c_{51}^2c_{55} -c_{51}^2c_{55}^2 -c_{41}c_{43}c_{52} -c_{45}c_{51}c_{52} 
+2c_{43}c_{51}c_{54} +c_{41}c_{52}c_{55} \\
&&+2c_{51}c_{54}c_{55} +c_{33} -c_{54}^2
\end{array}
$$
} 
\end{proposition}

\begin{proof}
Since we have a non-negative total arrow grading here, and since $I(\BO) \cap K[\Cnull]
= \{0\}$ by Proposition~\ref{prop-MaxDegBasis}.a, we can use~\cite[Alg.~4.6]{KR5}
to compute {\it all} best tuples of separating indeterminates.
We find 36 such tuples, and each of them has 13 elements. Using any of these tuples, 
we get a best $Z$-separating re-embedding of~$I(\BO)$.  
In other words, for each of them we get an isomorphism from~$B_\OO$ to a polynomial
ring in 12 indeterminates $C\setminus Z$ modulo an ideal generated by two polynomials. 
Among the 36 tuples, we selected one for which the corresponding pair of polynomials 
$(f_1,f_2)$ is such that $\#(\Supp(f_1)\cup \Supp(f_2))$ is minimal.
\end{proof}

Using the method of $Z$-separating re-embeddings, this is the best we can do.
However, by~\cite[Thm.~2.4]{Fog} and~\cite[Lemma 7]{Hui3}, we know that~$\BO$
is a regular scheme, and hence~\cite[Thm.~6.7]{KR5} implies that~$B_\OO$
is in fact isomorphic to a 10-dimensional polynomial ring.
Moreover, \cite[Alg.\ 6.8]{KR5} allows us to find an explicit isomorphism.
Let us apply it to the best $Z$-separating
re-embedding of the L-shape border basis scheme constructed above.

\begin{proposition}{\bf (Re-Embedding the L-Shape Border Basis Scheme)}\label{prop-theLshape}\\
Starting with the coordinate ring $B_\OO = K[C] / I(\BO)$ of the L-shape border basis scheme $\BO$, 
perform the following steps.
\begin{enumerate}
\item[(1)] Using Prop.~\ref{prop-L-Shape-Zsep}, construct a
re-embedding $\Phi: B_\OO \longrightarrow K[C\setminus Z]/ \langle f_1,\, f_2\rangle$.

\item[(2)] Define a $K$-algebra isomorphism $\psi:\; K[C \setminus Z] \longrightarrow
K[C\setminus Z]$ be letting 
\begin{align*}
\psi(c_{21}) &=\; -c_{41}^2c_{43}c_{51}^2 -c_{41}c_{45}c_{51}^3 +c_{41}^2c_{42}c_{51} +c_{41}^3c_{52} 
   +c_{41}^2c_{51}c_{54} -c_{45}c_{51}^2\\ 
&\quad  +c_{41}c_{51}c_{55} +c_{41}c_{42} -c_{41}c_{54} -c_{21} +c_{43}\\
\psi(c_{22}) &=\; 2c_{33}c_{41}c_{51} +2c_{35}c_{51}^2 +2c_{22} -c_{33}\\
\psi(c_{33}) &=\; c_{33}c_{41}c_{51} +c_{35}c_{51}^2 +c_{22}\\
\psi(c_{42}) &=\; c_{41}^2c_{43}c_{51}^3 +c_{41}c_{45}c_{51}^4 -c_{41}^2c_{42}c_{51}^2 -c_{41}^3c_{51}c_{52} 
   -c_{41}^2c_{51}^2c_{54} +c_{45}c_{51}^3\\ 
&\quad  -c_{41}c_{51}^2c_{55} -c_{41}c_{42}c_{51} +c_{41}c_{51}c_{54} +c_{21}c_{51} -c_{43}c_{51} -c_{42}
\end{align*}
and $\psi(c_{ij})=c_{ij}$ for the remaining indeterminates. This yields an induced $K$-algebra
isomorphism 
$$
\Psi:\; K[C\setminus Z]/ \langle f_1,\, f_2\rangle \;\longrightarrow\;
K[C \setminus Z] / \langle \psi(f_1),\, \psi(f_2)\rangle
$$
where $\psi(f_1)=c_{21}$ and $(\psi(f_1),\psi(f_2))$ is $(c_{21},c_{22})$-separating.

\item[(3)] Let $\hat{f}_1=c_{21}$ and substitute $c_{21} \mapsto 0$ in $\psi(f_2)$ to get
a polynomial $\hat{f}_2$. Then $\langle \psi(f_1), \psi(f_2)\rangle = \langle \hat{f}_1,
\hat{f}_2\rangle$ and the coherently $(c_{21},c_{22})$-separating tuple $(\hat{f}_1,\hat{f}_2)$
yields a separating re-embedding
$$
\Theta:\; K[C\setminus Z]/ \langle \hat{f}_1,\, \hat{f}_2\rangle \;\longrightarrow\;
K[\widehat{C}]
$$
where $\widehat{C} = (c_{33},c_{35},c_{41},c_{42},c_{43},c_{45},c_{51},c_{52},c_{54},c_{55})$.
\end{enumerate}
Altogether, the $K$-algebra isomorphism $\Theta\circ\Psi\circ\Phi:\; B_\OO \longrightarrow
K[\widehat{C}]$ implies that the L-shape border basis scheme is isomorphic to 
the affine space $\mathbb{A}^{10}$.
\end{proposition}

\begin{proof}
After the first step, we have to continue with the methods based on the unimodular matrix
problem we developed in~\cite[Sect.~6]{KR5}. More precisely, we have to apply~\cite[Alg.~6.8]{KR5}
which, in turn, is based on~\cite[Alg.~6.5]{KR5}. Let us follow the main steps.

The ring $K[\Cnull] = K[c_{41}, c_{51}]$ is the $K$-subalgebra of $K[C]$ 
generated by the indeterminates of total arrow degree zero. 
The $K[\Cnull]$-linear parts of~$f_1$ and~$f_2$ are
\begin{align*}
L_1 &= f_1 \;=\; (c_{41}^2c_{51}^2 + c_{41}c_{51} +1 )\, c_{21} +(c_{41}^2c_{51} +c_{41})\, c_{42}  
+(c_{41}^2c_{51}^2 -1 )\, c_{43}   \\
&  \qquad  +( c_{41}c_{51}^3 +c_{51}^2)\, c_{45}+(c_{41}^3)\, c_{52}  +(-c_{41}^2c_{51} +c_{41})\, c_{54}   
-(c_{41}c_{51} )\, c_{55} \\
L_2 &= (c_{41}c_{51})\,  c_{22}  +(-2c_{41}c_{51} +1)\, c_{33} -(c_{51}^2)\, c_{35}
\end{align*}
where $\deg_W(L_1)=1$ and $\deg_W(L_2)=2$.

In degree one, the coefficient matrix $A_1$ associated to $L_1$ is
{\small 
$$
A_1\!=\!( 
{-}c_{41}^2c_{51}^2 {-} c_{41}c_{51} {-}1 \ \  {-}c_{41}^2c_{51} {-} c_{41} \ \ 
{-}c_{41}^2c_{51}^2{+}1\ \  {-}c_{41}c_{51}^3 {-}c_{51}^2\ \ {-}c_{41}^3 \ \  
c_{41}^2c_{51}{-}c_{41} \ \ c_{41}c_{51}
)
$$
} 
and a solution of its UMP is 
$$
B_1 = 
\left( \begin{smallmatrix}
-1    &\ c_{41}^2c_{51} +c_{41}   & \  -c_{41}^2c_{51}^2 +1    &     \  -c_{41}c_{51}^3 -c_{51}^2        
         &\   c_{41}^3    &\  c_{41}^2c_{51} -c_{41} &\    {c_{41}c_{51}}_{\mathstrut}   \cr
c_{51}&\  -c_{41}^2c_{41}^2 -c_{41}c_{51} -1 &\  c_{41}^2c_{51}^3 -c_{51}       &\ c_{41}c_{51}^4 +c_{51}^3    
         &\ {-c_{41}^3c_{51}}^{\mathstrut}  &\  -c_{41}^2c_{51}^2 +c_{41}c_{51} &\ -c_{41}c_{51}^2 \cr
0&\  0&\      1& \      0&\      0&\        0&\    0 \cr
0&\  0&\      0& \      1&\      0&\        0&\    0 \cr
0&\  0&\      0& \      0&\      1&\        0&\    0\cr
0&\  0&\      0& \      0&\      0&\        1&\    0 \cr
0&\  0&\      0& \      0&\      0&\        0& \    1
\end{smallmatrix}\right)
$$
The multiplication of ~$B_1$ by $(c_{21}, c_{42}, c_{43}, c_{45}, c_{52}, c_{54}, c_{55})^{\rm tr}$ 
(whose entries are the indeterminates of $\deg_W= 1$) yields a column matrix $C_1$ of polynomials 
of total arrow degree one.

In total arrow degree two, the coefficient matrix $A_2$ of~$L_2$ is
$$
A_ 2 \;=\;   \left(
c_{41}c_{51}\;\; -2c_{41}c_{51}+1\;\; -c_{51}^2
\right)
$$
A solution of  the corresponding UMP is
$$
B_2 = \left( \begin{smallmatrix}
  2 \     & 2c_{41}c_{51}-1  \  & 2c_{51}^2  \cr
  1  \    &  c_{41} c_{51}    \   & c_{51}^2  \cr
  0  \   &       0                   \   &  1
\end{smallmatrix}\right)
$$
The multiplication of~$B_2$ by $(c_{22}, c_{33}, c_{35})^{\rm tr}$ (whose entries are the 
indeterminates of $\deg_W= 2$) yields a column matrix $C_2$ of polynomials of total arrow
degree two.

Integrating $C_1 \cup C_2$ into $C\setminus Z$ provides us with the $K$-algebra
automorphism~$\psi$. It is constructed such that $\psi(f_1),\psi(f_2))$
is $(c_{21},c_{22})$-separating. Thus we can make this pair coherently
separating in step~(3) and get the $(c_{21},c_{22})$-separating re-embedding~$\Theta$
which concludes the proof.
\end{proof}

In geometric jargon, we express the result of this proposition 
by saying that~$\BO$ is an {\it affine cell}.

\begin{remark}
The isomorphism 
$$
\Theta\circ\Psi\circ\Phi:\; K[C]/I(\BO) \;\longrightarrow\; 
K[c_{33},c_{35},c_{41},c_{42},c_{43},c_{45},c_{51},c_{52},c_{54},c_{55}]
$$ 
maps the residue classes of the indeterminates in~$C$
to polynomials whose supports have the lengths
$$
(78, 329, 375, 372, 419, 10, 87, 87, 95, 109, 8, 90, 86, 99,
1, 1, 11, 1, 11,  1, 1, 1, 9,  1, 1)
$$
Displaying them would exceed the limits of this paper.
The interested readers may find them in the \cocoa\ code available from the authors.
\end{remark}

\bigskip\bigbreak
%
%

\section{Simplicial Border Basis Schemes}
\label{sec5}

In this section we consider the following special kind of border basis scheme.
As usual, we let $P=K[x_1,\dots,x_n]$ be the polynomial ring over a field~$K$.

\begin{definition}\label{def-simplicial}
Let $\OO$ be an order ideal in~$\mathbb{T}^n$. Given $i\ge 0$, we write~$\OO_i$
for the set of terms of degree~$i$ in~$\OO$.
\begin{enumerate}
\item[(a)] Let $d\ge 1$. The order ideal $\OO$ is called {\bf simplicial of type~$d$}
if $\# \OO_i = \HF_P(i)$ for $i\le d$ and $\# \OO_i = 0$ for $i>d$.

\item[(b)] If~$\OO$ is a simplicial order ideal, the corresponding border
basis scheme~$\BO$ is called a {\bf simplicial border basis scheme}.

\end{enumerate}
\end{definition}

The following is a consequence of Proposition~\ref{prop-CharMaxDeg}.

\begin{corollary}\label{cor-CharSimplicial}
Let $\OO$ be an order ideal in~$\mathbb{T}^n$. Then the following conditions
are equivalent.
\begin{enumerate}
\item[(a)] The order ideal~$\OO$ is simplicial.

\item[(b)] The total arrow degree is a positive grading on~$K[C]$.

\end{enumerate}
\end{corollary}

\begin{proof}
Under both conditions, the order ideal~$\OO$ has a MaxDeg border by Proposition~\ref{prop-CharMaxDeg}. 
Condition~(b) means that $\deg(b_j)>d$ for all $j=1,\dots,\nu$. Since
Proposition~\ref{prop-CharMaxDeg} yields $\deg(t_i) \le d$ for $i=1,\dots,\mu$, 
condition~(b) is equivalent to the condition that all terms of degree~$d$ are in~$\OO$,
and hence to~(a).
\end{proof}

In view of the fact that the total arrow degree grading is positive
for simplicial order ideals, we know by~\cite[Cor.~4.8]{KR5}
that there exists an optimal separating re-embedding of~$I(\BO)$.
Our goal in this section is to describe an optimal separating re-embedding
of~$I(\BO)$ explicitly. More precisely, we show that a properly chosen
subset of the natural set of generators yields a separating tuple
which defines an optimal re-embedding of~$I(\BO)$.

As a small preparatory work, let us collect some basic numerical information 
about simplicial order ideals.

\begin{lemma}\label{lem-formulas}
Let $d\ge 1$, let $\OO = \{t_1, \dots, t_\mu\}$ be the simplicial order ideal of type~$d$
in~$\mathbb{T}^n$, and let $\partial\OO = \{ b_1, \dots, b_\nu\}$ be its border.
Then the following formulas hold.
\begin{enumerate}
\item[(a)] $\mu = \binom{d+n}{n}$

\item[(b)] $\nu = \binom{d+n}{n-1}$

\item[(c)] $\# \OOint = \binom{d+n-1}{n}$

\item[(d)] $\# \OOrim = \binom{d+n-1}{n-1}$

\item[(e)] $\# C = \binom{d+n}{n} \cdot \binom{d+n}{n-1}$

\item[(f)] $\# \Cint = \binom{d+n-1}{n} \cdot \binom{d+n}{n-1}$

\item[(g)] $\# \Crim = \binom{d+n-1}{n-1} \cdot \binom{d+n}{n-1}$
\end{enumerate}
\end{lemma}

\begin{proof}
These formulas follow straightforwardly from the 
affine Hilbert function of a polynomial ring given in~\cite[Props. 5.1.13 and 5.6.3.c]{KR2}.
\end{proof}

The next proposition implies that we can eliminate the indeterminates in~$\Cint$ 
using a $Z$-separating re-embedding and find a separating tuple among the natural 
generators of~$I(\BO)$.

\begin{proposition}\label{prop-simplicialZ}
Let $d\ge 1$, let $\OO= \{t_1, \dots, t_\mu\}$ be the simplicial order ideal of type~$d$
in~$\mathbb{T}^n$, and let $s = \#\Cint = \#\OOint \cdot \#(\partial\OO)=
\binom{d+n-1}{n} \cdot \binom{d+n}{n-1}$.
\begin{enumerate}
\item[(a)] There exists a term ordering $\sigma$ such that for each indeterminate $c_{ij} \in \Cint$ there
is a natural generator $f_{ij}$ of~$I(\BO)$ with $c_{ij}=\LT_\sigma(f_{ij})$.

\item[(b)] Using $Z=\Cint$, the tuple $F=(f_{ij})$ is $Z$-separating and consists of
natural generators of~$I(\BO)$.

\end{enumerate} 
\end{proposition}

\begin{proof}
First we prove claim~(a). Let $Z$ be the subtuple of~$C$ consisting of the 
indeterminates in~$\Cint$. 
We define an elimination ordering $\sigma$ for~$Z$ by constructing a suitable matrix~$M$
and letting $\sigma = {\tt Ord}(M)$ (see~\cite[Sect.~1.4]{KR1}).

In the first row of~$M$ and the column corresponding to an indeterminate 
$c_{ij}\in C$, we put the number $2\delta -1$, where~$\delta$ is the
total arrow degree $\deg_W(c_{ij}) = \deg(b_j) - \deg(t_i) =  d+1 -\deg(t_i)$.

Next we define the second row of the matrix $M$ as follows. For each arrow degree~$\delta$,
let $S_\delta$ be the set of all indeterminates $c_{ij}\in \Cint$ with
$\deg_W(c_{ij})=\delta$. Let $D$ be the set of all arrow degrees $\delta$ such that
$S_\delta\ne \emptyset$. For each $\delta = (\delta_1,\dots,  \delta_n)$ in~$D$, 
we choose a fixed,  for instance the first, number $k\in \{1,\dots,  n\}$ such that $\delta_k > 0$, 
and note  that this is possible by Proposition~\ref{prop-poscomp}. 
For $c_{ij} \in S_\delta$, we write
$t_i = x_1^{\alpha_1} \cdots x_n^{\alpha_n}$ with $\alpha_1,\dots, \alpha_n\in\mathbb{N}$.
Then we put the number $\alpha_k$ into the column corresponding to~$c_{ij}$
of the second row of the matrix~$M$. In the columns corresponding to $c_{ij}\in \Crim$
we put zeros into the second row of~$M$. Finally, we complete the matrix $M$ to
a matrix defining an elimination ordering $\sigma = {\tt Ord}(M)$ for~$Z$.

Our goal is to show that$\sigma$ is a $Z$-separating term ordering for~$I(\BO)$.
For this purpose we prove that for each $c_{ij} \in \Cint$ there exists a
generator $f_{ij}\in \AR_{\OO}$ such that the following conditions are satisfied:
\begin{enumerate} 
\item[(1)] The linear part of~$f_{ij}$ is (up to sign) $c_{ij}$ or of the form
$c_{ij} - c_{i'j'}$ with $c_{ij} >_\sigma c_{i'j'}$.

\item[(2)] For every term $c_{\alpha \beta} c_{\gamma\delta}$ of degree two in 
$\Supp(f_{ij})$, we have $c_{ij} >_\sigma c_{\alpha\beta} c_{\gamma\delta}$.
\end{enumerate}
Then it follows that $\LT_\sigma(f_{ij}) = c_{ij}$, and hence~$\sigma$ is a
$Z$-separating term ordering for the ideal $I(\BO)$.

To prove the existence of~$f_{ij}$ with the desired properties,
we use the explicit descriptions of the linear and quadratic parts
of the natural generators of~$I(\BO)$ provided in~\cite[Props.~6.3 and 6.6]{KLR3}. 
Let $c_{ij} \in \Cint$, and let $\delta = (\delta_1, \dots, \delta_n)$
be the arrow degree of~$c_{ij}$. Let $k\in \{1, \dots, n\}$ be the component of~$\delta$
with $\delta_k >0$ chosen in the above construction, and let $t_i = x_1^{\alpha_1}
\cdots x_n^{\alpha_n}$. Now we distinguish two cases:

{\it Case 1:} $\alpha_k=0$. By assumption, the exponent of~$b_j$ with respect to~$x_k$
is positive. Let $\ell\in \{1, \dots, n\} \setminus \{k\}$. Then $x_\ell \, b_j / x_k$ 
is an element of $\partial\OO$, and hence equal to $b_{j'}$ for some $j'\in \{1, \dots, \nu\}$.
Notice that $(b_j, b_{j'})$ is an across-the-rim neighbor pair.
Moreover, the term $t_m = x_\ell t_i$ is contained in~$\OO$, since $t_i\in \OOint$.
However, the term~$t_m$ is not divisible by~$x_k$. Next, consider the polynomial
$f_{ij} = \AR(j,j')_m$. The linear part of~$f_{ij}$ is (up to sign) equal to~$c_{ij}$. 
Furthermore, the terms in the quadratic part of~$f_{ij}$ are of the form 
$c_{m\kappa} c_{\rho j''}$ where $j''\in \{j,j'\}$
and $c_{\rho j''} \in \Crim$. Here $\deg(t_m)= \deg(t_i)+1$ implies that the entry 
in the first row of~$M$ corresponding to $c_{m\kappa}$ is two less than the entry 
corresponding to~$c_{ij''}$.
Since the entry corresponding to $c_{\rho j''}$ is one, it follows that we have
$c_{ij} >_\sigma c_{m\kappa} c_{\rho j''}$. Thus we get $\LT_\sigma(f_{ij}) = c_{ij}$.

{\it Case 2:} $\alpha_k >0$. In this case, we choose an index $\ell \in \{1, \dots, n\}
\setminus \{k\}$ and note that $x_\ell b_j / x_k$ is in~$\partial\OO$. Let $j'\in\{1, \dots, \nu\}$
be such that $x_k b_{j'} = x_\ell b_j$. Since $t_i \in \OOint$, we know that $t_m = x_\ell t_i$
is in~$\OO$. As we are in the case $\alpha_k>0$, it follows that~$x_k$ divides~$t_m$, 
and there exists an index $i'\in \{1, \dots, \mu\}$ with $x_k t_{i'} = t_m = x_\ell t_i$.
Then $(b_j, b_{j'})$ is an across-the-rim neighbor pair, and we consider the polynomial
$f_{ij} = \AR(j,j')_m$. The linear part of~$f_{ij}$ is (up to sign) given by
$c_{ij} -c_{i' j'}$. Here the facts that the total arrow degree of~$c_{ij}$ and~$c_{i' j'}$
are the same and that the $x_k$-exponent of~$t_{i'}$ is $\alpha_k -1$ imply
$c_{ij} >_\sigma c_{i'j'}$. 

The terms of degree two in the support of~$f_{ij}$ are seen to be smaller than~$c_{ij}$
with respect to~$\sigma$ in exactly the same way as in Case~1. Altogether, it follows
that $\LT_\sigma(f_{ij}) = c_{ij}$, and the proof of~(a) is complete.

Claim~(b) follows immediately from the proof of~(a).
\end{proof}

This proposition has several useful consequences.

\begin{corollary}
Let $d\ge 1$, and let $\OO = \{t_1, \dots, t_\mu\}$ be the simplicial order ideal of type~$d$
in~$\mathbb{T}^n$. We denote the maximal ideal of~$K[C]$ generated by the indeterminates
by $\M = \langle c_{ij}\rangle$ and its residue class ideal in~$B_\OO$ 
by $\m = \langle \bar{c}_{ij} \rangle$.
\begin{enumerate}
\item[(a)] The linear part of~$I(\BO)$ is given by
$\Lin_\M(I(\BO)) = \langle \Cint \rangle_K$.

\item[(b)] We have $\dim_K (\Cot_\m(B_\OO))= \# \Crim = \binom{d+n-1}{n-1} 
\cdot \binom{d+n}{n-1}$. 

\item[(c)] We have $\edim(B_\OO) = \dim_K (\Cot_{\m}(B_\OO)) $.  
\end{enumerate}
In particular, for the term ordering~$\sigma$ constructed in the proposition, 
the $Z$-sep\-a\-ra\-ting re-embedding  
$\Phi: K[C] / I(\BO) \longrightarrow K[\Crim] / ( I\cap K[\Crim])$
defined by reducing the $\sigma$-Gr\"obner basis~$F$ is optimal.
\end{corollary}

\begin{proof}
To prove(a), we note that the indeterminates of~$\Cint$ are in
$\Lin_\M(I(\BO))$ by the proposition. By~\cite[Cor.~2.8.b]{KSL}, there is no
linear part of one of the natural generators of~$I(\BO)$ which contains
an indeterminate of~$\Crim$ in its support. Hence we get
$\Lin_\M(I(\BO)) = \langle \Cint\rangle_K$.

The first equality in~(b) follows directly from the proposition and~(a).
Since we have the equality $\#C = \#\OO \cdot \#(\partial\OO) =
\binom{d+n}{n}\cdot \binom{d+n}{n-1}$, we get
$$
\# \Crim = \#C - \#\Cint =   \tbinom{d+n}{n}\cdot \tbinom{d+n}{n-1} 
- \tbinom{d+n-1}{n} \cdot \tbinom{d+n}{n-1}
= \tbinom{d+n-1}{n-1}\cdot \tbinom{d+n}{n-1}
$$
Finally, we prove~(c). As an application of~\cite[Cor. 4.2]{KLR1}, we know that
the embedding dimension of~$B_\OO$ is $\#C - \#\Cint = \#\Crim$, 
and this number is equal to $\dim_K (\Cot_\m(B_\OO))$ by~(b).
\end{proof}

Another application of the proposition is that we can check easily whether 
a simplicial border basis scheme is an affine cell.

\begin{proposition}\label{prop-AffineOrSing}
Let $d\ge 1$, let $\OO$ be the simplicial order ideal of type~$d$, and
let $p=(0, \dots, 0)$ be the monomial point of~$\BO$.
\begin{enumerate}
\item [(a)] If $n=2$ then $\BO$ is isomorphic to $\mathbb{A}^{(d+1)(d+2)}$.

\item[(b)] If $n\ge 3$ then $p$ is a singular point of~$\BO$.
\end{enumerate}
\end{proposition}

\begin{proof}
For $n = 2$, it follows from Proposition~\ref{prop-simplicialZ} that
$\BO$ can be embedded into an affine space of dimension
$\dim_K (\Cot_\m(B_\OO)) = \binom{d+1}{1}\cdot \binom{d+2}{1} =   (d+1)(d+2)$.
It is also known that the  scheme $\BO$ is irreducible and its dimension is given by
the formula
$\dim(\BO)=2\cdot  \#(\OO) = 2\frac{(d+1)(d+2)}{2} = (d+1)(d+2)$.
The conclusion follows.

To prove(b), we notice that the dimension of the principal component of~$\BO$
is $n\cdot \#\OO = n\cdot \binom{d+n}{n}$.
Moreover, from Proposition~\ref{prop-simplicialZ}.c we get
$\dim_K (\Cot_\m(B_\OO)) = \binom{d+n-1}{n-1} \cdot \binom{d+n}{n-1}$.
Hence it suffices to show that for $n\ge 3$ we have
$$
n\cdot \tbinom{d+n}{n} < \tbinom{d+n-1}{n-1} \cdot \tbinom{d+n}{n-1}  \eqno{(1)}
$$
To prove (1), we note that $n\cdot \binom{d+n}{n} {=}(d+n)\cdot \binom{d+n-1}{n-1}$.
Hence(1) is equivalent to
$$
d+n <  \tbinom{d+n}{n-1}   \eqno{(2)}
$$
Therefore we need to show that for $n\ge 3$ we have
$$
1 <  \tfrac{ (d+n-1)\cdot(d+n-2)\cdots(d+2)}{(n-1)!}=
\tfrac{d+n-1}{n-1}\cdot \tfrac{d+n-2}{n-2}\cdots \tfrac{d+2}{2}
$$
which is obvious, since $d\ge 1$.

Notice that for $n=2$ we have $d+n = \binom{d+n}{n-1}$. Thus we have equality in(2),
and hence also equality in~(1). This shows $\dim_K(\Cot_\m(B_\OO)) = \dim(\BO)$ 
and implies that~$p$ is a smooth point of~$\BO$, in agreement with~(a).
\end{proof}

The following examples illustrate the main results of this section.

\begin{example}\label{ex-123}
In $P=\mathbb{Q}[x,y]$, consider the simplicial order ideal~$\OO$
of type~2, i.e., $\OO = \{1,y,x,y^2,xy,x^2\}$. Its border is 
$\partial\OO = \{y^3, xy^2, x^2y,x^3\}$ (see Figure~\ref{fig:123simp}).

\begin{figure}[ht]
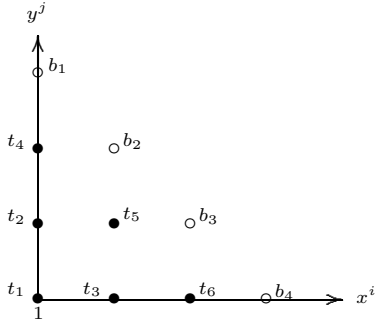

\centering{\makebox{
\beginpicture
		\setcoordinatesystem units <1cm,1cm>
		\setplotarea x from 0 to 4, y from 0 to 3.5
		\axis left /
		\axis bottom /
		\arrow <2mm> [.2,.67] from  3.5 0  to 4 0
		\arrow <2mm> [.2,.67] from  0 3  to 0 3.5
		\put {$\scriptstyle x^i$} [lt] <0.5mm,0.8mm> at 4.1 0.1
		\put {$\scriptstyle y^j$} [rb] <1.7mm,0.7mm> at 0 3.6
		\put {$\bullet$} at 0 0
		\put {$\bullet$} at 1 0
		\put {$\bullet$} at 0 1
		\put {$\bullet$} at 0 2
		\put {$\bullet$} at 1 1
		\put {$\bullet$} at 2 0

		\put {$\scriptstyle 1$} [lt] <-1mm,-1mm> at 0 0
		\put {$\scriptstyle t_1$} [rb] <-1.3mm,0.4mm> at 0 0
		\put {$\scriptstyle t_3$} [rb] <-1.3mm,0.4mm> at 1 0
		\put {$\scriptstyle t_2$} [rb] <-1.3mm,0mm>   at 0 1
		\put {$\scriptstyle t_4$} [rb] <-1.3mm,0mm>   at 0 2
		\put {$\scriptstyle t_5$} [lb] <0.8mm,0.4mm>  at 1 1
		\put {$\scriptstyle t_6$} [lb] <0.8mm,0.4mm>  at 2 0
		\put {$\scriptstyle b_1$} [rb] <4.3mm,0mm> at 0 3
		\put {$\scriptstyle b_2$} [lb] <0.8mm,0mm> at 1 2
		\put {$\scriptstyle b_3$} [lb] <0.8mm,0mm> at 2 1
		\put {$\scriptstyle b_4$} [lb] <0.8mm,0mm> at 3 0
		
		\put {$\circ$} at 0 3
		\put {$\circ$} at 2 1
		\put {$\circ$} at 1 2
		\put {$\circ$} at 3 0
\endpicture
}} 
\caption{A simplicial order ideal and its border}\label{fig:123simp}
\end{figure}

Then the simplicial border basis 
scheme~$\BO$ is naturally embedded in $\mathbb{A}^{24}$.
Using Lemma~\ref{lem-formulas}, we calculate 
$\# \OOint = \binom{d+n-1}{n} = 3$,
$\# \OOrim = \binom{d+n-1}{n-1} = 3$,
$\# \Cint = \binom{d+n-1}{n} \cdot \binom{d+n}{n-1} = 12$, and
$\# \Crim = \binom{d+n-1}{n-1} \cdot \binom{d+n}{n-1} = 12$.
Hence we let $Z=\Cint$ and apply the $Z$-separating re-embedding of 
Proposition~\ref{prop-simplicialZ}.b to get $\BO \cong \mathbb{A}^{12}$.
\end{example}

In the final example of this section we construct an optimal embedding of a border basis scheme
which is not an isomorphism with an affine space.

\begin{example}\label{ex-13}
In $P=\QQ[x,y,z]$, consider the simplicial order ideal $\OO = \{1,z,y,x\}$ of type~1.
Here we have $\# \OOint = \binom{d+n-1}{n} = 1$, $\# \OOrim = \binom{d+n-1}{n-1} = 3$,
$\# \Cint = \binom{d+n-1}{n} \cdot \binom{d+n}{n-1} = 6$, and
$\# \Crim = \binom{d+n-1}{n-1} \cdot \binom{d+n}{n-1} = 18$.

Therefore Proposition~\ref{prop-ZsepReemb} yields a term ordering~$\sigma$ and 
a tuple of natural generators~$F$ of~$I(\BO)$ such that $Z=\Cint$ is exactly the 
tuple of $\sigma$-leading terms of~$F$. 
The reduced $\sigma$-Gr\"obner basis of $\langle F \rangle$ 
defines an optimal re-embedding of~$B_\OO$. 
More specifically, it defines an isomorphism between $B_\OO$ and a ring of the form
$\QQ[C\setminus Z]/J$, where $\#(C\setminus Z)=18$ and where the ideal~$J$ is
minimally generated by 15 homogeneous quadratic polynomials (see~\cite[Example~5.6]{KLR1}).

Thus we have found an optimal re-embedding of the 12-dimensional simplicial border basis 
scheme~$\BO$ into $\mathbb{A}^{18}$.
As predicted by Proposition~\ref{prop-AffineOrSing}.b, the monomial point $p=(0,\dots,0)$ 
is a singular point of~$\BO$.
\end{example}

\bigskip\bigbreak
%
%

\section{Planar Border Basis Schemes}
\label{sec6}

Given an order ideal~$\OO$ in two indeterminates, we say
that the corresponding scheme $\BO$ is a {\bf planar border basis scheme}.
Recall that we usually write $x=x_1$ and $y=x_2$ in this case, so we
are working with a polynomial ring $P=K[x,y]$ over a field~$K$.

In this section we examine general planar border basis schemes,
and in the final two sections of the paper we look at some particular cases.
Recall that planar border basis schemes are always regular by~\cite[Thm.~2.4]{Fog}
and~\cite[Lemma 7]{Hui3}.
Our first goal is to find a $Z$-separating re-embedding, where~$Z$ is chosen as large
as possible in general. We need the following auxiliary terminology.

\begin{definition} 
Let $\OO = \{ t_1,\dots,  t_\mu \}$ be an order ideal in~$\mathbb{T}^2$ with border
$\partial\OO = \{b_1,\dots, b_\nu \}$.
\begin{enumerate}
\item[(a)] A term is called {\bf $x$-free} (resp. {\bf $y$-free}) if~$x$ (resp.~$y$)
does not divide it.

\item[(b)] A term $b_j\in\partial\OO$ is said
to have an {\bf up-neighbor} $b_{j'}\in\partial\OO$ if
$b_{j'} = x b_j$ or $b_{j'} = y b_j$. In this case~$b_j$
is also called a {\bf down-neighbor} of~$b_{j'}$.

\item[(c)] A sequence of border terms $b_{j_1},\dots, b_{j_k}$ 
is said to form a {\bf plateau} if the pairs
$(b_{j_1},b_{j_2})$, $\dots$, $(b_{j_{k-1}}, b_{j_k})$ form a maximal
chain of across-the-rim neighbor pairs with $x b_{j_\ell} = y b_{j_{\ell+1}}$
for $\ell = 1,\dots, k-1$, and if $b_{j_1}, b_{j_k}$ have no
up-neighbors.

\item[(d)] Let $b_{j_1},\dots, b_{j_k}$ be a plateau, 
and let $b_{j_1} = x b_{j'_1}$ with $b_{j'_1}\in\partial\OO$.
Let $b_{j'_2},\dots, b_{j'_m}$ be a maximal sequence of border terms
such that $(b_{j'_{\ell+1}}, b_{j'_\ell})$ are across-the-rim or next-door neighbors
in $x$-direction, i.e., such that $x b_{j'_{\ell+1}} \in 
\{ b_{j'_\ell}, y b_{j'_\ell} \}$. Then the sequence $b_{j'_1},\dots, 
b_{j'_m}$ is called the {\bf $x$-leg} of the  plateau $b_{j_1},\dots, b_{j_k}$.
The {\bf $y$-leg} of the plateau $b_{j_1},\dots, b_{j_k}$ is defined
analogously.

\end{enumerate}
\end{definition}

Let us illustrate this definition in an easy case.

\begin{example}\label{ex-122-2} 
In $P=\mathbb{Q}[x,y]$, look at the order ideal $\OO = \{t_1,\dots, t_5\}$, where 
$t_1=1$, $t_2=y$, $t_3=x$, $t_4=y^2$, and $t_5=xy$. Then the border of~$\OO$
is $\partial\OO = \{b_1,\dots, b_4\}$, where $b_1=x^2$, $b_2=y^3$, 
$b_3=xy^2$, and $b_4=x^2y$.

\begin{figure}[ht]
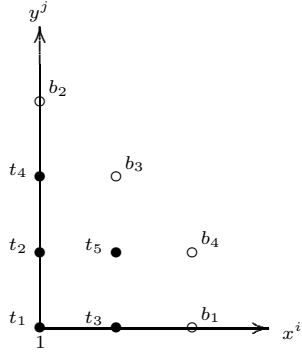

\centering{\makebox{
\beginpicture
		\setcoordinatesystem units <1cm,1cm>
		\setplotarea x from 0 to 3, y from 0 to 3.5
		\axis left /
		\axis bottom /
		\arrow <2mm> [.2,.67] from  2.5 0  to 3 0
		\arrow <2mm> [.2,.67] from  0 3.5  to 0 4
		\put {$\scriptstyle x^i$} [lt] <0.5mm,0.8mm> at 3.1 0
		\put {$\scriptstyle y^j$} [rb] <1.7mm,0.7mm> at 0 4
		\put {$\bullet$} at 0 0
		\put {$\bullet$} at 1 0
		\put {$\bullet$} at 0 1
		\put {$\bullet$} at 1 1
		\put {$\bullet$} at 0 2
		\put {$\scriptstyle 1$} [lt] <-1mm,-1mm> at 0 0
		\put {$\scriptstyle t_1$} [rb] <-1.3mm,0.4mm> at 0 0
		\put {$\scriptstyle t_3$} [rb] <-1.3mm,0.4mm> at 1 0
		\put {$\scriptstyle t_2$} [rb] <-1.3mm,0mm> at 0 1
		\put {$\scriptstyle t_4$} [rb] <-1.3mm,0mm> at 0 2
		\put {$\scriptstyle t_5$} [rb] <-1.3mm,0mm> at 1 1
		\put {$\scriptstyle b_1$} [lb] <0.8mm,0.8mm> at 2 0
		\put {$\scriptstyle b_2$} [lb] <0.8mm,0.8mm> at  0 3
		\put {$\scriptstyle b_3$} [lb] <0.8mm,0.8mm> at 1 2
	    \put {$\scriptstyle b_4$} [lb] <0.8mm,0.8mm> at 2 1
	
		\put {$\circ$} at 2 0
		\put {$\circ$} at 2 1
		\put {$\circ$} at 1 2
		\put {$\circ$} at 0 3
\endpicture
}} 
\caption{An order ideal with one plateau}\label{fig:122-2}
\end{figure}

There is only one border term with an up-neighbor, namely $b_1=x^2$.
The terms $(b_2,b_3,b_4)$ form a plateau which has
the $y$-leg $(b_1)$ and no $x$-leg.
\end{example}

For planar border basis schemes, our next goal is to reprove and make explicit 
some beautiful results from~\cite{Hui1}. 
More precisely, let~$Z$ be the tuple of non-exposed 
indeterminates $C \setminus \Cexp$. To show that
there exists a $Z$-separating re-embedding of~$I(\BO)$ which uses the
natural set of generators of~$I(\BO)$, we proceed as follows:
\begin{enumerate}
\item[(1)] Define a set of weights for all $c_{ij}$ such that
$\wt(c_{ij}) = 0$ if $c_{ij}$ is an exposed indeterminate, such that
$\wt(c_{ij}) > 0$ if $c_{ij}$ is a non-exposed indeterminate,
and such that for each non-exposed indeterminate $c_{ij}$ there
exists a natural generator~$f$ of~$I(\BO)$ for which $c_{ij}$
is the unique term of highest weight in~$\Supp(f)$.

\item[(2)] Define an elimination ordering~$\sigma$ for~$Z$
given by a matrix whose first row consists of the weights assigned
in Step~(1).
\end{enumerate}

Then, by~\cite[Prop.~4.2 and Def.~4.4]{KLR2}, we have a $Z$-separating
re-em\-bed\-ding of~$I(\BO)$. In this way, we have eliminated all non-exposed
indeterminates. To execute Step~(1), we exhibit the following explicit algorithm.
For its formulation, we use the definitions of the tuples $\ND(j,j')$ 
and $\AR(j,j')$ from \cite[Def.~6.1]{KLR3} and the relation $\sim$ to denote
cotangent equivalence (see~\cite[Sect.~5]{KLR3}).

\begin{algorithm}{\bf (Weight Assignment Algorithm)}\label{alg-WeightAssign}\\
Let $\mathcal{O} \subseteq \mathbb{T}^2$ be an order ideal in the indeterminates~$x,y$,
let $I(\BO) \subseteq K[C]$ be the ideal defining the $\mathcal{O}$-border basis scheme, 
and let $\delta = \max \{ \deg(t_i) \mid t_i\in \OO\}$.
Consider the following instructions.
\begin{enumerate}
\item[(1)] Let $S=\emptyset$, and let $\wt(c_{ij}) = 0$ for all $c_{ij}\in \Cexp$.

\item[(2)] For $d=\delta, \delta-1, \dots, 0$, execute the loop defined by 
steps (3)--(7).

\item[(3)] Let $C_d$ be the set of all indeterminates~$c_{ij}$ in~$C$
such that $\deg(t_i)=d$ and such that $c_{ij}$ has not yet been assigned
a weight.

\item[(4)] Let $c_{ij} \in C_d$ be an indeterminate for which~$b_j$
has an up-neighbor $b_{j'} = x_k b_j$.
Let $i'\in\{1,\dots, \mu\}$ such that $t_{i'} = x_k t_i$, and let
$b_{j_1},\dots, b_{j_m}$ be the $x_k$-exposed border terms.
Then we set
$$
\wt(c_{ij}) \;=\; 1 + \wt(c_{i'j'}) + \wt(c_{i' j_1}) + \dots + \wt(c_{i' j_m})
$$

\item[(5)] Let $c_{ij} \in C_d$ be an indeterminate for which~$b_j$
is an element of a plateau $b_{j_1},\dots, b_{j_k}$. To define
the weight of~$c_{ij}$, we distinguish several cases.
In steps (6a) to (6e), we treat the cases where $\deg_W(c_{ij})$
has a positive $x$-component. In steps (7a) to (7e), we treat the cases
where $\deg_W(c_{ij})$ has a non-positive $x$-component and a positive 
$y$-component.

\item[(6a)] Suppose that $j=j_\ell$ for some $\ell\in\{2,\dots, k\}$
and $t_i$ is $x$-free. Moreover, let $m\in\{1,\dots, \mu\} $ be such that
$t_m = y t_i$. Then we set
$$
\wt(c_{ij}) \;=\; 1 + \wt(c_{m1}) + \cdots + \wt(c_{m\nu})
$$

\item[(6b)] Suppose there exists $\ell\in\{2,\dots, k\}$
such that $j=j_\ell$ and such that we have $t_i = x^\alpha y^\beta$
with $\alpha\le \ell -2$. Starting with $\ell=2$, the case treated in
step~(6a), we define $\wt(c_{ij})$ by induction on~$\ell$.
Let $i'\in\{1,\dots, \mu\}$ such that $c_{i' j_{\ell-1}} \sim c_{i j_\ell}= c_{ij}$
and let $m\in\{1,\dots, \mu\}$ such that $t_m = y t_i$. Then we define
$$
\wt( c_{ij} ) \;=\; 1 + \wt(c_{i' j_{\ell-1}}) + \wt(c_{m1}) + \cdots + \wt(c_{m\nu})
$$

\item[(6c)] Suppose that $j=j_1$ and that the plateau has an $x$-leg at~$b_{j_1}$
starting with $b_{j_1} = x b_{j'_1}$. Moreover, suppose that $t_i$ is $x$-free.
Let $b_{j''_1},\dots, b_{j''_m}$ be the $x$-exposed border terms.
Then we introduce a new indeterminate $p_{ij}$, we set $\wt(c_{ij}) = p_{ij}$
and we append the following inequalities to~$S$: 
$p_{ij} > \wt(c_{i j''_1})$, $\dots$, $p_{ij} > \wt(c_{i j''_m})$.

\item[(6d)] [Going down the $x$-leg]\\
Suppose that $j=j_1$ and that $b_{j'_1},\dots, b_{j'_\kappa}$ is the $x$-leg
of the plateau. Then we introduce a new indeterminate $p_{ij}$ and
we set $\wt(c_{ij}) = p_{ij}$. Let $\ell\le \kappa$ and let 
$i'_1,\dots, i'_\ell \in\{1,\dots, \mu\}$ be those indices for which 
$c_{i'_1 j'_1} \sim \cdots \sim c_{i'_\ell j'_\ell}$. Then we set
$$
\wt(c_{i'_1 j'_1}) \;=\; p_{ij}-1,\;\; \wt(c_{i'_2 j'_2}) \;=\; p_{ij}-2,\;\dots, \;
\wt(c_{i'_\ell j'_\ell}) \;=\; p_{ij}-\ell
$$
and we add inequalities of the form $p_{ij} \ge N$ to~$S$ 
to express the fact that $c_{i'_\lambda j'_\lambda}$ is the term of the 
highest weight in the polynomial $\ND(j'_{\lambda+1}, j'_\lambda)_{i'_\lambda}$
resp. $\AR(j'_\lambda, j'_{\lambda+1})_{m_\lambda}$.

\item[(6e)] Suppose that $j=j_\ell$ for some $\ell\in\{ 1,\dots, k\}$
and suppose that $t_i = x^\alpha y^\beta$ with $\alpha \ge \ell-1$.
Starting with the case $\ell=1$, treated in
step~(6c), we define $\wt(c_{ij})$ by induction on~$\ell$.
Choose $i'\in\{1,\dots, \mu\}$ such that we have $c_{i' j_{\ell-1}} \sim c_{i j_\ell}= c_{ij}$,
and let $m\in\{1,\dots, \mu\}$ such that $t_m = y t_i$. Then we define
$$
\wt( c_{ij} ) \;=\; 1 + \wt(c_{i' j_{\ell-1}}) + \wt(c_{m1}) + \cdots + \wt(c_{m\nu})
$$

\item[(7a)] Suppose that $j=j_\ell$ for some $\ell\in\{1,\dots, k-1\}$
and $t_i$ is $y$-free. Let $m\in\{1,\dots, \mu\}$ be such that
$t_m = x t_i$. Then we set
$$
\wt(c_{ij}) \;=\; 1 + \wt(c_{m1}) + \cdots + \wt(c_{m\nu})
$$

\item[(7b)] Suppose there exists $\ell\in\{1,\dots, k-1\}$
such that $j=j_\ell$ and such that $t_i = x^\alpha y^\beta$
with $\beta\le k-1-\ell$. Starting with the case $\ell=k-1$, treated in
step~(7a), we define $\wt(c_{ij})$ by downward induction on~$\ell$.
Moreover, let $i'\in\{1,\dots, \mu\}$ be such that $c_{i' j_{\ell+1}} \sim c_{i j_\ell}= c_{ij}$,
and let $m\in\{1,\dots, \mu\}$ be such that $t_m = x t_i$. Then we define
$$
\wt( c_{ij} ) \;=\; 1 + \wt(c_{i' j_{\ell+1}}) + \wt(c_{m1}) + \cdots + \wt(c_{m\nu})
$$

\item[(7c)] Suppose that $j=j_k$ and that the plateau has a $y$-leg at~$b_{j_k}$
starting with $b_{j_k} = y b_{j'_1}$. Moreover, suppose that $t_i$ is $y$-free.
Let $b_{j''_1},\dots, b_{j''_m}$ be the $y$-exposed border terms.
Then we introduce a new indeterminate $p_{ij}$, we set $\wt(c_{ij}) = p_{ij}$
and append the following inequalities to~$S$: 
$p_{ij} > \wt(c_{i j''_1})$, $\dots$, $p_{ij} > \wt(c_{i j''_m})$.

\item[(7d)] [Going down the $y$-leg]\\
Suppose that $j=j_k$ and that $b_{j'_1},\dots, b_{j'_\kappa}$ is the $y$-leg
of the plateau. Then we introduce a new indeterminate $p_{ij}$ and
we set $\wt(c_{ij}) = p_{ij}$. Let $\ell\le \kappa$ and let 
$i'_1,\dots, i'_\ell \in\{1,\dots, \mu\}$ be those indices for which 
$c_{i'_1 j'_1} \sim \cdots \sim c_{i'_\ell j'_\ell}$. Then we set
$$
\wt(c_{i'_1 j'_1}) \;=\; p_{ij}-1,\;\; \wt(c_{i'_2 j'_2}) \;=\; p_{ij}-2,\;\dots, \;
\wt(c_{i'_\ell j'_\ell}) \;=\; p_{ij}-\ell
$$
and we add inequalities of the form $p_{ij} \ge N$ to~$S$ 
to express the fact that $c_{i'_\lambda j'_\lambda}$ is the term of the 
highest weight in the polynomial $\ND(j'_\lambda, j'_{\lambda+1})_{i'_\lambda}$
resp. $\AR(j'_\lambda, j'_{\lambda+1})_{m_\lambda}$.

\item[(7e)] Suppose that $j=j_\ell$ for some $\ell\in\{ 1,\dots, k\}$
and suppose that $t_i = x^\alpha y^\beta$ with $\beta \ge k-\ell$.
Starting with the case $\ell=k$ which was treated in
step~(7c), we define $\wt(c_{ij})$ by induction on~$\ell$.
Choose $i'\in\{1,\dots, \mu\}$ such that $c_{i' j_{\ell+1}} \sim c_{i j_\ell}= c_{ij}$,
and let $m\in\{1,\dots, \mu\}$ such that $t_m = x t_i$. Then we define
$$
\wt( c_{ij} ) \;=\; 1 + \wt(c_{i' j_{\ell+1}}) + \wt(c_{m1}) + \cdots + \wt(c_{m\nu})
$$

\item[(8)] Finally, consider the set of inequalities~$S$. Choose the numbers $p_{ij}$
large enough so that all inequalities are satisfied. In this way, all 
indeterminates~$c_{ij}$ will be assigned explicit weights.
\end{enumerate}
This is an algorithm which assigns to each $c_{ij}\in C$
a weight $\wt(c_{ij})\ge 0$ such that the following conditions are satisfied:
\begin{enumerate}
\item[(a)] All exposed indeterminates $c_{ij}\in \Crim$ satisfy
$\wt(c_{ij}) = 0$.

\item[(b)] All non-exposed indeterminates $c_{ij} \in \Cint$ satisfy
$\wt(c_{ij}) > 0$.

\item[(c)] For each non-exposed indeterminate $c_{ij}$ there
exists a natural generator~$f$ of the ideal~$I(\BO)$ for which $c_{ij}$
is the unique term of highest weight in~$\Supp(f)$.
\end{enumerate}
\end{algorithm} 

\begin{proof}
Since the algorithm is clearly finite, it suffices to show that
the assigned weights lead to the desired terms of highest weight
in specifically chosen natural generators of~$I(\BO)$.
\begin{enumerate}
\item[(1)] The definition of the weights in step~(4) implies that $c_{ij}$
is the highest weight term in $\ND(j,j')_{i'}$.

\item[(2)] In case (6a), the term $c_{ij} = c_{i j_\ell}$ is the linear part
and the highest weight term of the polynomial $\AR(j_\ell,j_{\ell-1})_m$.

\item[(3)] In case (6b), the linear part of $\AR(j_\ell,j_{\ell-1})_m$ is
(up to sign) the polynomial $c_{ij} - c_{i' j_{\ell-1}}$ and the highest
weight term is $c_{ij}$.

\item[(4)] In case (6c), the term $c_{ij}$ is the linear part of 
$\ND(j'_1,j_1)_i$ and its weight~$p_{ij}$ is the highest weight of
any term in this polynomial by the given inequalities for~$p_{ij}$.

\item[(5)] In case (6d), when we move from~$b_{j'_\lambda}$ to
$b_{j'_{\lambda+1}}$ while going down the $x$-leg, we go down a next-door
neighbor pair using the polynomial $\ND(j'_{\lambda+1}, j'_\lambda)_{i'_\lambda}$
or an across-the-rim neighbor pair using the polynomial 
$\AR(j'_\lambda,j'_{\lambda+1})_{m_\lambda}$ and weight given 
to~$c_{i_\lambda j_\lambda}$ as well as the given inequalities for $p_{ij}$
make sure that the term of highest weight is $c_{i_\lambda j_\lambda}$
in both cases.

\item[(6)] In case (6e), we are again using the linear part 
of $\AR(j_\ell,j_{\ell-1})_m$ to move along the plateau. It is
(up to sign) the polynomial $c_{ij} - c_{i' j_{\ell-1}}$ and the highest
weight term is~$c_{ij}$ by the given definition.

\item[(7)] In steps (7a)-(7e) analogous arguments apply.
\end{enumerate}
Finally, we note that the required set of inequalities for the
weights $p_{ij}$ can always be solved: when we move into and down an $x$-leg in
steps (6c) and (6d), the weight $p_{ij}$ has to be at least as high as 
the weights $\wt(c_{i'j'})$ with $i'=i$ or $b_{j'}$ in the $x$-tail of the
plateau. On the other hand, the inequalities in steps (7c) and (7d) involve
indices~$i$ which correspond to the other end of the plateau or its $y$-leg.
Therefore we will never get inequalities which interfere with each other,
and by choosing $p_{ij}$ large enough all inequalities will be satisfied.
\end{proof}

When we compare this algorithm to the 
ones in~\cite[Alg.~3.2 and 4.1]{AKL}, we notice that it is restricted to the planar case and to
the tuple~$Z$ of non-exposed indeterminates.
However, it has the advantage of always returning a result,
whereas the more general algorithms may fail to come up with suitable weights.

Based on the weight assignment algorithm and our strategy, we have now shown the following
general result.

\begin{theorem}{\bf (Elimination of Non-Exposed Indeterminates)}\label{thm-elimNex}\\
Let $K$ be a field, let $P=K[x,y]$,
let $\OO = \{ t_1,\dots, t_\mu \}$ be an order ideal in~$\mathbb{T}^2$, 
let $\partial\OO = \{b_1,\dots, b_\nu \}$ be its border, and let
$C = \{ c_{ij} \mid i=1,\dots, \mu;\; j=1,\dots, \nu \}$  
denote the set of coefficients of the generic $\OO$-border prebasis.
Let $Z=(z_1,\dots, z_s)$ be the tuple consisting of the non-exposed 
indeterminates $C\setminus \Cexp$.

Using Algorithm~\ref{alg-WeightAssign}, determine
weights $p_{ij}\in\mathbb{N}$ for all $c_{ij}\in C$, and let~$\sigma$
be an elimination ordering for~$Z$ defined by putting the weights~$p_{ij}$
into the first row of a square matrix of size $\mu\nu$ and extending it to
a term ordering matrix.

Then the natural generators $f_1,\dots, f_s$ of~$I(\BO)$ given in the proof of 
Algorithm~\ref{alg-WeightAssign} satisfy $\LT_\sigma(f_i) = z_i$
for $i=1,\dots, s$ and define a $Z$-separating re-embedding of~$I(\BO)$.
In other words, the non-exposed indeterminates have been eliminated. 
\end{theorem}

\begin{proof}
This follows from the algorithm and~\cite[Prop.~2.2]{KLR2}.
\end{proof}

This theorem has the following immediate consequences.
Recall that an indeterminate $c_{ij}$ is called {\bf basic} if it is not
contained in the linear part of a generator of~$I(\BO)$ (see~\cite[Def.~5.1 and Lemma 5.2.b]{KLR3}).

\begin{corollary}\label{cor-exposedgenerate}
In the setting of the theorem, the following claims hold.
\begin{enumerate}
\item[(a)] The residue classes modulo $I(\BO)$ of the exposed indeterminates generate 
the $K$-algebra $B_\OO$.

\item[(b)] Every basic indeterminate is exposed and cannot be contained in a tuple~$Z$
for which $I(\BO)$ has a $Z$-separating re-embedding.

\end{enumerate}
\end{corollary}

\begin{proof} 
Claim (a) follows immediately from the theorem, and claim~(b) 
is a consequence of~(a) and~\cite[Thm.~5.6.a]{KLR3}.
\end{proof}

Let us compare Theorem~\ref{thm-elimNex} and its corollary to~\cite[Prop.~7.1.2]{Hui1}.

\begin{remark}\label{rem-Hui1approach}
In~\cite[Sect.~6 and~7]{Hui1}, the author uses an ingenious argument
based on the Hilbert-Burch Theorem to show that the residue classes modulo $I(\BO)$
of the exposed indeterminates form a $K$-algebra system of generators of~$B_\OO$.
In particular, using the maximal minors of the matrix~$M$ whose rows
represent the liftings of the ND and AR neighbor syzygies to 
$\Syz_{B_\OO}(\bar{g}_1,\dots,  \bar{g}_\nu)$, where $G=\{ g_1,\dots, g_\nu\}$
is the generic $\OO$-border prebasis, 
one can find explicit representations of the non-exposed indeterminates~$c_{ij}$ 
as polynomial expressions $c_{ij} = f_{ij}(\Cexp)$ in the exposed indeterminates. 
The advantage of this method is that those expressions form
a Gr\"obner basis of the ideal $\langle c_{ij} - f_{ij}(\Cexp)\rangle$,
and this Gr\"obner basis is obtained by expanding the minors of~$M$.

On the other hand, the advantage of the approach based on the 
Weight Assignment Algorithm is that it finds a subset of the natural set of 
generators of~$I(\BO)$ which yields the re-embedding, and these tend to be rather 
sparse polynomials. However, it is clear that both approaches pay the price that 
the polynomials in the corresponding reduced Gr\"obner basis, which is necessary to 
compute the re-embedding, may have large supports.
\end{remark}

Of particular interest are order ideals $\OO$ such that the non-exposed indeterminates
produce an isomorphism of~$B_\OO$ with a polynomial ring. They can be described
as follows.

\begin{remark}\label{rem-optZsepPlanar}
In the above setting, suppose that the $Z$-separating re-embedding
of~$I(\BO)$ is optimal. By~\cite[Cor.~2.9.b]{KLR3}, this implies that 
the border basis scheme~$\BO$ is isomorphic to $\mathbb{A}^{2\mu}_K$.
Moreover, we have $\# Z  = \# C - 2\mu = \mu \nu - 2\mu$.

The special cases in which the $Z$-separating re-embedding is optimal for 
the tuple~$Z$ of non-exposed indeterminates, i.e., the cases where $\# \Cexp = 2\mu$, 
are classified in~\cite[Prop.~7.3.1]{Hui1}.
They occur if and only if the order ideal~$\OO$ has the ``sawtooth'' form described there.
However, as already mentioned in~\cite[Rem.~7.5.3.1]{Hui1}, there are cases 
where~$\BO$ is an affine cell, but $\OO$ does not have the ``sawtooth'' form (see 
for instance~\cite[Example~6.8]{KLR3}).
\end{remark}

Let us see an example which illustrates this remark.

\begin{example}{\bf (The (2,3)-Box)}\label{ex-23Exposed}\\
In $P=\mathbb{Q}[x,y]$, consider the order ideal 
$\OO = \{t_1,t_2,t_3,t_4, t_5, t_6\}$, where 
$t_1=1$, $t_2=y$, $t_3=x$, $t_4=y^2$, $t_5=xy$, and $t_6=xy^2$. Its border
is $\partial\OO = \{b_1, b_2, b_3, b_4, b_5\}$, where $b_1=x^2$, $b_2=y^3$, $b_3 = x^2y$,
$b_4 = xy^3$, and $b_5 = x^2y^2$.  
\begin{figure}[ht]
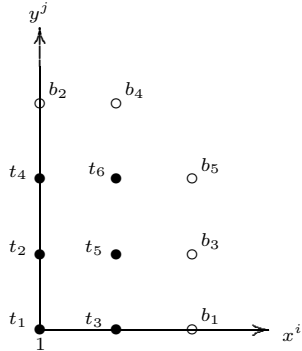

\centering{\makebox{
\beginpicture
		\setcoordinatesystem units <1cm,1cm>
		\setplotarea x from 0 to 3, y from 0 to 3.5
		\axis left /
		\axis bottom /
		\arrow <2mm> [.2,.67] from  2.5 0  to 3 0
		\arrow <2mm> [.2,.67] from  0 3.5  to 0 4
		\put {$\scriptstyle x^i$} [lt] <0.5mm,0.8mm> at 3.1 0
		\put {$\scriptstyle y^j$} [rb] <1.7mm,0.7mm> at 0 4
		\put {$\bullet$} at 0 0
		\put {$\bullet$} at 1 0
		\put {$\bullet$} at 0 1
		\put {$\bullet$} at 1 1
		\put {$\bullet$} at 1 2
		\put {$\bullet$} at 0 2
		\put {$\scriptstyle 1$} [lt] <-1mm,-1mm> at 0 0
		\put {$\scriptstyle t_1$} [rb] <-1.3mm,0.4mm> at 0 0
		\put {$\scriptstyle t_3$} [rb] <-1.3mm,0.4mm> at 1 0
		\put {$\scriptstyle t_2$} [rb] <-1.3mm,0mm> at 0 1
		\put {$\scriptstyle t_4$} [rb] <-1.3mm,0mm> at 0 2
		\put {$\scriptstyle t_5$} [rb] <-1.3mm,0mm> at 1 1
		\put {$\scriptstyle t_6$} [lb] <-4mm,0.4mm> at 1 2
		\put {$\scriptstyle b_1$} [lb] <0.8mm,0.8mm> at 2 0
		\put {$\scriptstyle b_2$} [lb] <0.8mm,0.8mm> at  0 3
		\put {$\scriptstyle b_3$} [lb] <0.8mm,0.8mm> at 2 1
	    \put {$\scriptstyle b_4$} [lb] <0.8mm,0.8mm> at 1 3
	    \put {$\scriptstyle b_5$} [lb] <0.8mm,0.8mm> at 2 2
	
		\put {$\circ$} at 2 0
		\put {$\circ$} at 2 1
		\put {$\circ$} at 2 2
		\put {$\circ$} at 1 3
		\put {$\circ$} at 0 3
\endpicture
}} 
\caption{The $(2,3)$-box order ideal and its border}\label{fig:2-3-box}
\end{figure}

The elements of the upper border are~$b_2$ and~$b_4$, and the elements of the upper rim
are~$t_4$ and~$t_6$. The elements of the right border are $b_1$, $b_3$, and~$b_5$. 
The elements of the right rim are~$t_3$, $t_5$, and~$t_6$.

The next-door neighbor pair $x b_2 = b_4$ yields the $x$-exposed indeterminates 
$c_{32}$, $c_{52}$, and $c_{62}$. The across-the-street neighbor pair $x b_4 = y b_5$ 
yields the $x$-exposed indeterminates $c_{34}$, $c_{54}$, and~$c_{64}$.
Hence the set of $x$-exposed indeterminates is $\{ c_{32}, c_{34}, c_{52}, c_{54},
c_{62}, c_{64} \}$. Analogously, the set of $y$-exposed indeterminates is 
$\{ c_{41}, c_{43}, c_{45}, c_{61}, c_{63}, c_{65}\}$. Altogether, we have found 12 exposed
indeterminates and we conclude that there is an isomorphism 
$$
B_\OO \cong \QQ[c_{32}, c_{34}, c_{41}, c_{43}, c_{45}, c_{52}, c_{54}, 
c_{61}, c_{62}, c_{63}, c_{64}, c_{65}]
$$
The weights provided by Algorithm~\ref{alg-WeightAssign} are
$$
\begin{array}{ccccccccccccccc}
c_{11}&  c_{12}&  c_{13}&  c_{14}&  c_{15}&  c_{21}&  c_{22}&  c_{23}&  c_{24}&  c_{25}&  
   c_{31}&  c_{32}&  c_{33}&  c_{34}&  c_{35}  \cr
13& 15& 13& 20& 19& 3& 5& 3& 4& 3& 2& 0& 3& 0& 9 \cr
c_{41}&  c_{42}&  c_{43}&  c_{44}&  c_{45}&  c_{51}&  c_{52}&  c_{53}&  c_{54}&  c_{55}&  
   c_{61}&  c_{62}&  c_{63}&  c_{64}&  c_{65} \cr
  0&1& 0& 1& 0& 1& 0& 1& 0& 2& 0& 0& 0& 0& 0
\end{array}
$$
and the algorithm computes the 18 natural generators whose leading terms are the indeterminates 
with positive weight, i.e., the non-exposed indeterminates.
\end{example}

Let us point out a special merit of Algorithm~\ref{alg-WeightAssign}.

\begin{remark}\label{rem:goodweights}
Not all orderings which eliminate non-exposed indeterminates find suitable natural generators
having such indeterminates as leading terms.
For instance, in the above example the polynomial $f\!=\! -c_{22}c_{41} -c_{24}c_{61} -c_{11}+c_{23}$
is a natural generator of the ideal~$I(\BO)$. The weights of the indeterminates 
appearing in~$\Supp(f)$ are 
$$
\begin{matrix}
c_{11}&    c_{22}&  c_{23}&  c_{24}& c_{41}&  c_{61}\cr
   13&  5& 3& 4& 0 & 0
\end{matrix}
$$ 
Hence we have $\LT_\sigma(f) = c_{11}$. Instead, if we use an elimination ordering~$\tau$ 
represented by a matrix whose first row assigns~1 to the non-exposed indeterminates
and 0 to the exposed indeterminates, we get $\LT_\tau(f) = c_{22}c_{41}$. 
Therefore~$\tau$ does not find 18 natural generators having the correct leading terms.
\end{remark}

Next we make good use of exposed indeterminates in the process of 
searching for optimal re-embeddings.  This goal is achieved by a suitable variant 
of~\cite[Alg.~5.7]{KLR3}. It uses the following terminology (cf.~\cite{KSL} and~\cite{KLR3}).

\begin{definition}\label{def-cotequiv}
Let $\m = \langle c_{ij} + I(\BO) \mid i\in \{1,\dots,\mu\},\, j\in \{1,\dots,\nu\} \rangle$
be the irrelevant maximal ideal of $B_\OO = K[C]/I(\BO)$. 
For every indeterminate $c_{ij}\in C$, let $\bar{c}_{ij}$ denote
its residue class in the cotangent space $\m / \m^2$ 
of~$B_\OO$ at the origin.
\begin{enumerate}
\item[(a)] The equivalence relation~$\sim$ on~$C$ defined
by $c_{ij} \sim c_{k\ell} \Leftrightarrow \langle \bar{c}_{ij}\rangle_K = 
\langle\bar{c}_{k\ell}\rangle_K$ is called {\bf cotangent equivalence}.

\item[(b)] The {\bf trivial} cotangent equivalence class consists of
all indeterminates $c_{ij}$ such that $\bar{c}_{ij}=0$.

\item[(c)] If a cotangent equivalence class contains at least two distinct elements, 
it is called {\bf proper}.

\end{enumerate}
\end{definition}

Now we adapt~\cite[Alg.~5.7]{KLR3} to get a version which computes
for every planar border basis scheme all optimal $Z$-separating
re-embeddings, i.e., the $Z$-separating re-embeddings which reach the embedding dimension.
Notice that the following algorithm may not find an optimal $Z$-separating re-embedding of~$\BO$ at all.

\begin{algorithm}{\bf (Optimal $Z$-Separating Re-embeddings)}\label{alg-OptPlanarZ}\\
Let $\OO = \{ t_1, \dots, t_\mu\}$ be an order ideal in~$\mathbb{T}^2$ with border
$\partial \OO = \{ b_1,\dots, b_\nu\}$, and let $C=(c_{ij})$ be the tuple of indeterminates 
which are the coefficients of the generic $\OO$-border prebasis. 
Consider the following sequence of instructions.
\begin{enumerate}
\item[(1)] Compute the trivial cotangent equivalence class~$E_0$
and also the proper cotangent equivalence classes $E_1,\dots,E_q$.

\item[(2)] Let $S=\emptyset$, and let $\Cexp$ be the tuple of exposed indeterminates in~$C$.
Moreover, for $i=1, \dots, q$, compute $\widetilde{E}_i =  E_i \cap \Cexp$.

\item[(3)] From each set $E_0\cup (C\setminus \Cexp) \cup \widetilde{E}_1^\ast \cup \cdots 
\cup \widetilde{E}_q^\ast$, where $\widetilde{E}_i^\ast$ is obtained from~$\widetilde{E}_i$ 
by deleting one element, form a tuple~$Z$ and perform the following steps. 

\item[(4)] Check whether the
ideal~$I(\BO)$ is $Z$-separating. If it is, append~$Z$ to~$S$.

\item[(5)]  Continue with Step~(3) and form the next tuple~$Z$ until 
all tuples have been dealt with. Then return~$S$ and stop.

\end{enumerate}
Then the result is an algorithm which computes the set~$S$ of all tuples~$Z$ 
of distinct indeterminates in~$C$ such that the
$Z$-separating re-embedding of~$I(\BO)$ is optimal.
\end{algorithm}

\begin{proof}
From~\cite[Alg.~5.7]{KLR3}, we know that the indeterminates 
in~$E_0$ must be eliminated. From Theorem~\ref{thm-elimNex}, we know that non-exposed 
indeterminates can be eliminated. Therefore, if we want to eliminate them, 
we do not need to use the sets~$E_i$, as we do in~\cite[Alg.~5.7]{KLR3}.
Instead, it suffices to use the sets~$\widetilde{E}_i$. Hence the correctness of the algorithm
follows from the correctness of~\cite[Alg.~5.7]{KLR3}.
\end{proof}

To perform step~(4) of this algorithm, we can use one of the methods of~\cite[Rem.~2.5]{KLR3}
or try the one of the algorithms~\cite[Alg.~3.2 and~4.1]{AKL}.
The following example, already mentioned in \cite[Example 6.8]{KLR3},
shows it at work.

\begin{example}\label{ex-1232}
In $P =\QQ[x,y]$, consider the order ideal
$\OO = \{t_1,\dots,t_8\}$ given by $t_1=1$, $t_2=y$, $t_3=x$, $t_4=y^2$, $t_5=xy$, $t_6=x^2$, 
$t_7=y^3$, and $t_8 =xy^2$. Its border is $\partial\OO = \{b_1,\dots,b_5\}$, where $b_1=x^2y$, 
$b_2=x^3$, $b_3=y^4$, $b_4 = xy^3$, and $b_5=x^2y^2$. 
\begin{figure}[ht]
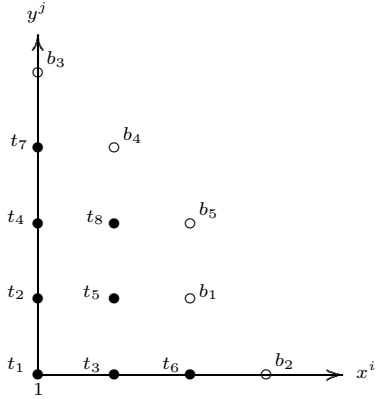

\centering{\makebox{
\beginpicture
		\setcoordinatesystem units <1cm,1cm>
		\setplotarea x from 0 to 4, y from 0 to 4.5
		\axis left /
		\axis bottom /
		\arrow <2mm> [.2,.67] from  3.5 0  to 4 0
		\arrow <2mm> [.2,.67] from  0 4  to 0 4.5
		\put {$\scriptstyle x^i$} [lt] <0.5mm,0.8mm> at 4.1 0.1
		\put {$\scriptstyle y^j$} [rb] <1.7mm,0.7mm> at 0 4.6
		\put {$\bullet$} at 0 0
		\put {$\bullet$} at 1 0
		\put {$\bullet$} at 0 1
		\put {$\bullet$} at 1 1
	    \put {$\bullet$} at 1 2
		\put {$\bullet$} at 0 2
		\put {$\bullet$} at 0 3
		\put {$\bullet$} at 2 0

		\put {$\scriptstyle 1$} [lt] <-1mm,-1mm> at 0 0
		\put {$\scriptstyle t_1$} [rb] <-1.3mm,0.4mm> at 0 0
		\put {$\scriptstyle t_3$} [rb] <-1.3mm,0.4mm> at 1 0
		\put {$\scriptstyle t_2$} [rb] <-1.3mm,0mm> at 0 1
		\put {$\scriptstyle t_4$} [rb] <-1.3mm,0mm> at 0 2
		\put {$\scriptstyle t_5$} [rb] <-1.3mm,0mm> at 1 1
		\put {$\scriptstyle t_6$} [lb] <-4mm,0.4mm> at 2 0
		\put {$\scriptstyle t_7$} [lb] <-4mm,0mm> at 0 3
		\put {$\scriptstyle t_8$} [lb] <-4mm,0mm> at 1 2
		\put {$\scriptstyle b_1$} [lb] <0.8mm,0mm> at 2 1
		\put {$\scriptstyle b_2$} [lb] <0.8mm,0.8mm> at 3 0
	    \put {$\scriptstyle b_3$} [lb] <0.8mm,0.8mm> at 0 4
	    \put {$\scriptstyle b_4$} [lb] <0.8mm,0.8mm> at 1 3
	    \put {$\scriptstyle b_5$} [lb] <0.8mm,0.8mm> at  2 2

		\put {$\circ$} at 3 0
		\put {$\circ$} at 2 2
		\put {$\circ$} at 2 1
		\put {$\circ$} at 1 3
		\put {$\circ$} at 0 4
\endpicture
}} 
\caption{A non-sawtooth order ideal and its border}\label{fig:non-sawtooth}
\end{figure}

Thus $\QQ[C] =\QQ[c_{11}, \dots, c_{85}]$ is a polynomial ring in~$40$ indeterminates. 
Notice that the scheme~$\BO$ is smooth and has dimension $\dim(\BO) = 2\mu = 16$.
In~\cite[Example 6.8]{KLR3}, we calculated the cotangent equivalence classes
\begin{align*}
E_0 = & \{ c_{11},  c_{12},  c_{13},  c_{14},  c_{15},  c_{21},  c_{22},  c_{23},  
c_{24},  c_{25},  c_{31},  c_{32},  c_{33},  c_{34},  c_{35}, \\ 
& \; c_{42},  c_{44},  c_{45},  c_{55},  c_{65} \}\\
E_1 = & \{c_{51},  c_{85}\}, \qquad E_2 = \{c_{43},  c_{54}\}, \qquad
E_3 = \{c_{41},  c_{52},  c_{75}\}
\end{align*}
The exposed indeterminates are
$\Cexp = ( c_{51},  c_{53},  c_{54},  c_{61},  c_{62},  c_{63},  c_{64},  
c_{65}, c_{71},  c_{72},  \allowbreak   c_{73},  c_{74},  c_{75},  
c_{81},  c_{82},  c_{83},  c_{84},  c_{85} )$.

Notice that $c_{52}$ and~$c_{55}$ are rim indeterminates which are not exposed.
Moreover, as~$c_{65}$ is exposed and trivial, it will be eliminated. 
At this point we already know that $\Cexp$ does not generate $B_\OO$ minimally.
Furthermore, the indeterminates $c_{41}$, $c_{43}$, and $c_{52}$ are non-exposed. 
Therefore we have $\widetilde{E}_1 = E_1$,  $\widetilde{E}_2 =\{c_{54} \}$, 
and $\widetilde{E}_3 = \{c_{75} \}$, so that $\widetilde{E}_2^\ast =\emptyset$ 
and $\widetilde{E}_3^\ast =\emptyset$.
Consequently, the remaining tuples to be examined are only the two tuples~$Z_1$ and~$Z_2$
formed from $E_0 \cup (C \setminus \Cexp) \cup \{ c_{51}\}$ and 
$E_0 \cup (C\setminus \Cexp) \cup \{ c_{85}\}$, respectively.

It turns out that both provide optimal $Z$-separating re-embeddings
$\BO \cong \mathbb{A}^{16}$.
In the first case, the exposed indeterminates $c_{65}$ and $c_{51}$ are eliminated, 
and in the second case, the exposed indeterminates $c_{65}$ and $c_{85}$ are eliminated. 

As already noticed in Remark~\ref{rem-optZsepPlanar},
this is an example of an order ideal~$\OO$ which does not have the ``sawtooth'' form.
Nevertheless, we get an optimal re-embedding of $I(\BO)$, and $\BO$ turns out to be an affine cell.
Notice also that the example mentioned in~\cite[Rem.~ 7.5.3]{Hui1}, is similar to this example.  
\end{example}

\bigskip\bigbreak
%
%

\section{Special Planar Border Basis Schemes}
\label{sec7}

To conclude this paper, we consider two types of special planar border basis schemes:
box border basis schemes and MaxDeg border basis schemes. Let us begin by looking
at box border basis schemes.

\begin{definition}
Let $K$ be a field, let $P=K[x_1,\dots, x_n]$, and let $a_1,\dots, a_n\in\NN_+$.
\begin{enumerate}
\item[(a)] The order ideal 
$$
\OO = \mathbb{T}^n \setminus \langle x_1^{a_1}, \dots,  x_n^{a_n} \rangle =
\{ x_1^{i_1} \cdots x_n^{i_n} \mid 0\le i_j < a_j\;\hbox{for}\;
j=1,\dots, n\} 
$$
is called the {\bf box} of type $(a_1,\dots, a_n)$.

\item[(b)] If~$\OO$ is a box, the corresponding border basis
scheme is called a {\bf box border basis scheme}.
\end{enumerate}
\end{definition}

Using the Weight Assignment Algorithm~\ref{alg-WeightAssign}, we now show
that planar box border basis schemes are affine cells. The following
terminology will prove useful.

\begin{definition}
Let $a,b\in\mathbb{N}_+$, and let $\OO = \{ x^i y^j \mid
0\le i< a,\; 0\le j< b\}$ be the box of type $(a,b)$ in
$P=K[x,y]$.
\begin{enumerate}
\item[(a)] The set $\{ x^{a-1}y^k \mid 0\le k\le b-1\}$ is called the {\bf right rim}
of~$\OO$.

\item[(b)] The set $\{ x^k y^{b-1} \mid 0\le k\le a-1\}$ is called the {\bf upper rim}
of~$\OO$.

\item[(c)] The set $\{ x^a y^k \mid 0\le k\le b-1\}$ is called the {\bf right border}
of~$\OO$

\item[(d)] The set $\{x^k y^b \mid 0\le k\le a-1\}$ is called the {\bf upper border}
of~$\OO$.

\end{enumerate}
\end{definition}

At this point we are ready to prove the promised result.

\begin{proposition}{\bf (Re-Embedding Planar Box Border Basis Schemes)}\label{prop-PlanarBox}\\
Let $a,b\in\mathbb{N}_+$, and let $\OO = \{ x^i y^j \mid
0\le i< a,\; 0\le j< b\}$ be the box of type $(a,b)$ in
$P=K[x,y]$.
\begin{enumerate}
\item[(a)] An indeterminate~$c_{ij}$ is $x$-exposed if and only if
$t_i$ is a term in the right rim and $b_j$ is a term in the
upper border of~$\OO$.

In particular, the number of $x$-exposed indeterminates~$c_{ij}$ is~$ab$.

\item[(b)] An indeterminate $c_{ij}$ is $y$-exposed if and only if
$t_i$ is a term in the upper rim and~$b_j$ is a term in the
right border of~$\OO$.

In particular, the number of $y$-exposed indeterminates~$c_{ij}$ is~$ab$.

\item[(c)] The tuple~$Z$ of all non-exposed indeterminates is $I(\BO)$-separating.
The $Z$-sep\-a\-rating re-embedding yields an isomorphism
$\Phi:\; B_\OO \longrightarrow K[C\setminus \Cexp]$.

In particular, every planar box border basis scheme is an affine cell.

\end{enumerate}
\end{proposition}

\begin{proof}
To prove~(a), we consider the up-neighbors in $x$-direction, i.e.,
the next-door neighbor pairs $x b_j = b_{j'}$ for which $b_j, b_{j'}$
are in the upper border.
All terms $b_j$ in the upper border, except for $b_j=x^{a-1} y^b$,
have an up-neighbor in $x$-direction, i.e., there is a next-door neighbor pair
$x b_j = b_{j'}$. The term $b_j=x^{a-1} y^b$ has an across-the-street neighbor
$b_{j'} =  x^a y^{b-1}$, since $x b_j = y b_{j'}$.
Altogether, for each term~$b_j$ in the upper border and for each term~$t_i$ in 
the right rim we get an $x$-exposed indeterminate~$c_{ij}$. In this way 
we get all $x$-exposed indeterminates.

The proof of~(b) follows analogously. It remains to prove~(c).
Since there are no other exposed indeterminates, we have $\# \Cexp = 2ab$.
By the Weight Assignment Algorithm~\ref{alg-WeightAssign}, we know that we can 
eliminate all non-exposed indeterminates, and the exposed
indeterminates generate~$B_\OO$ (see Corollary~\ref{cor-exposedgenerate}.a). Therefore
$\dim(\BO) = 2 \mu = 2ab$ shows that $\BO \cong \mathbb{A}^{2ab}_K$,
and that~$B_\OO$ is a polynomial ring over~$K$ generated by the residue classes
of the indeterminates in $\Cexp$ described in~(a) and~(b).
\end{proof}

Since we have seen already several examples for this proposition (see
Examples~\ref{ex-box22}, \ref{ex-12exposed}, and~\ref{ex-23Exposed}), 
we refrain from presenting further ones.

\medskip

Our second topic in this section are planar MaxDeg border basis schemes.
By Corollary~\ref{cor-bivariateaffinecells}, we know already that they 
are affine cells if~$K$ is a perfect field. The necessary re-embedding may involve
an isomorphism which is constructed using a solution of the Unimodular Matrix Problem (see for 
instance Proposition~\ref{prop-theLshape}).
However, in the simplicial case we were able to construct the re-embedding directly
without assuming anything about the base field and by applying a $Z$-separating re-embedding only.

In the following we look at more cases when we can identify a planar MaxDeg border basis scheme
with an affine space with the aid of a $Z$-separating re-embedding.
From here on we assume that $\OO=\{t_1,\dots,t_\mu\}$ is a non-simplicial MaxDeg order ideal,
and we let $d=\max\{\deg(t_i) \mid i=1,\dots,\mu\}$.
Recall from Proposition~\ref{prop-CharMaxDeg} that~$\OO$ has a generic Hilbert function.
In our setting this implies
$\# \OO_i = \dim_K(P_i) = i+1$ for $i=0,\dots,d-1$ and
$1\le \#\OO_d \le d$, as well as $\#\OO_i = 0$ for $i>d$.

\begin{definition}\label{def-OdSegment}
Let $\OO$ be a non-simplicial MaxDeg order ideal in~$\mathbb{T}^2$, and let
$d = \max\{\deg(t) \mid t\in \OO\}$.
\begin{enumerate}
\item[(a)] Let $0\le i \le j\le d$. The tuple of terms $(x^i y^{d-i},\, x^{i+1} y^{d-i-1},\,
\dots, x^j y^{d-j})$ in $\OO_d^{j-i+1}$ is called an {\bf $\OO_d$-segment} if $i=0$ or $j=d$, or if
$x^{i-1} y^{d-i+1} \notin\OO$ and $x^{j+1} y^{d-j-1}\notin \OO$. 

\item[(b)] For an $\OO_d$-segment as in~(a), the number $\ell = j-i+1$ is called the
{\bf length} of the $\OO_d$-segment.

\item[(c)] The number $s\ge 1$ of $\OO_d$-segments of~$\OO$ is called the 
{\bf segmentation type} of~$\OO$.
\end{enumerate}
\end{definition}

Let us illustrate these concepts in a simple case.

\begin{example}\label{ex-somegaps}
Consider the following two non-simplicial MaxDeg order ideals in~$\mathbb{T}^2$
consisting of $\mu=8$ terms and having maximal degree $d=3$.
\begin{enumerate}
\item[(a)] The order ideal $\OO = \{ 1, x, y, x^2, xy, y^2, y^3, xy^2 \}$ has segmentation 
type~1. Its only $\OO_3$-segment is $(y^3, xy^2)$ and has length~2.

\item[(b)] The order ideal $\OO' = \{ 1, x, y, x^2, xy, y^2, y^3, x^2y \}$ has segmentation 
type~2. Its two $\OO'_3$-segments are $(y^3)$ and $(x^2y)$, and they are both of length~1.
\end{enumerate}

\begin{figure}[ht]
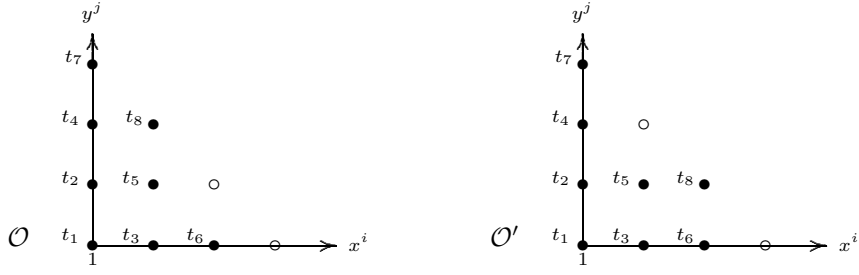

$\OO$\quad
\makebox{\beginpicture
		\setcoordinatesystem units <0.8cm,0.8cm>
		\setplotarea x from 0 to 4, y from 0 to 3.5
		\axis left /
		\axis bottom /
		\arrow <2mm> [.2,.67] from  3.5 0  to 4 0
		\arrow <2mm> [.2,.67] from  0 3  to 0 3.5
		\put {$\scriptstyle x^i$} [lt] <0.5mm,0.8mm> at 4.1 0.1
		\put {$\scriptstyle y^j$} [rb] <1.7mm,0.7mm> at 0 3.6
		\put {$\bullet$} at 0 0
		\put {$\bullet$} at 1 0
		\put {$\bullet$} at 0 1
		\put {$\bullet$} at 1 1
	    \put {$\bullet$} at 0 2
		\put {$\bullet$} at 2 0
		\put {$\bullet$} at 0 3
		\put {$\bullet$} at 1 2

		\put {$\scriptstyle 1$} [lt] <-1mm,-1mm> at 0 0
		\put {$\scriptstyle t_1$} [rb] <-1.3mm,0.4mm> at 0 0
		\put {$\scriptstyle t_3$} [rb] <-1.3mm,0.4mm> at 1 0
		\put {$\scriptstyle t_2$} [rb] <-1.3mm,0mm> at 0 1
		\put {$\scriptstyle t_4$} [rb] <-1.3mm,0mm> at 0 2
		\put {$\scriptstyle t_5$} [rb] <-1.3mm,0mm> at 1 1
		\put {$\scriptstyle t_6$} [lb] <-4mm,0.4mm> at 2 0
		\put {$\scriptstyle t_7$} [lb] <-4mm,0mm> at 0 3
		\put {$\scriptstyle t_8$} [lb] <-4mm,0mm> at 1 2
		
		\put {$\circ$} at 2 1
		\put {$\circ$} at 3 0
\endpicture
} 
\qquad\qquad
$\OO'$\quad
\makebox{\beginpicture
		\setcoordinatesystem units <0.8cm,0.8cm>
		\setplotarea x from 0 to 4, y from 0 to 3.5
		\axis left /
		\axis bottom /
		\arrow <2mm> [.2,.67] from  3.5 0  to 4 0
		\arrow <2mm> [.2,.67] from  0 3  to 0 3.5
		\put {$\scriptstyle x^i$} [lt] <0.5mm,0.8mm> at 4.1 0.1
		\put {$\scriptstyle y^j$} [rb] <1.7mm,0.7mm> at 0 3.6
		\put {$\bullet$} at 0 0
		\put {$\bullet$} at 1 0
		\put {$\bullet$} at 0 1
		\put {$\bullet$} at 2 0
	    \put {$\bullet$} at 1 1
		\put {$\bullet$} at 0 2
		\put {$\bullet$} at 0 3
		\put {$\bullet$} at 2 1

		\put {$\scriptstyle 1$} [lt] <-1mm,-1mm> at 0 0
		\put {$\scriptstyle t_1$} [rb] <-1.3mm,0.4mm> at 0 0
		\put {$\scriptstyle t_3$} [rb] <-1.3mm,0.4mm> at 1 0
		\put {$\scriptstyle t_2$} [rb] <-1.3mm,0mm> at 0 1
		\put {$\scriptstyle t_4$} [rb] <-1.3mm,0mm> at 0 2
		\put {$\scriptstyle t_5$} [rb] <-1.3mm,0mm> at 1 1
		\put {$\scriptstyle t_6$} [lb] <-4mm,0.4mm> at 2 0
		\put {$\scriptstyle t_7$} [lb] <-4mm,0mm> at 0 3
		\put {$\scriptstyle t_8$} [lb] <-4mm,0mm> at 2 1

        \put {$\circ$} at 1 2
        \put {$\circ$} at 3 0
\endpicture
} 
\caption{Two order ideals and their segments}\label{fig:segments}
\end{figure}

\end{example}

For non-simplicial MaxDeg order ideals in~$\mathbb{T}^2$, the number of indeterminates~$c_{ij}$
can be calculated as follows.

\begin{lemma}\label{lem-NuFormula}
Let $\OO \subseteq \mathbb{T}^2$ be a non-simplicial MaxDeg order ideal,
and let $d= \max\{ \deg(t) \mid t\in\OO \}$.
\begin{enumerate}
\item[(a)] Denoting the lengths of the $\OO_d$-segments by $\ell_1,\dots,\ell_s$,
we obtain the formula $\mu = d(d+1) + \ell_1 + \cdots + \ell_s$.

\item[(b)] The number~$\nu$ of elements of the border $\partial\OO$ is given
by $\nu=d+1+s$.
\end{enumerate}
In particular, the number of indeterminates in~$C=(c_{ij})$ is
$$
\mu\nu \;=\; [ d(d+1) + \ell_1+\cdots +\ell_s] \cdot (d+1+s)
$$
\end{lemma}

\begin{proof}
Claim (a) follows from $\#\OO = \#\OO_0 + \cdots + \#\OO_d = 
1 +\cdots + d + \ell_1 + \cdots +\ell_s$.

It remains to prove~(b). Clearly, there are $(d+1) - \ell_1 - \ell_s$
border terms of~$\OO$ in degree~$d$. It remains to show that
every $\OO_d$-segment $(x^i y^{d-i},\, x^{i+1} y^{d-i-1},\,
\dots,\, x^j y^{d-j})$ with $0\le i\le j \le d$ gives rise to
$\ell_s+1 = j-i+2$ border terms in degree~$d+1$. 
These are exactly the terms $b'_1 = x^i y^{d-i+1}$, $b'_2 = 
x^{i+1} y^{d-i}$, $\dots$, $b'_{j-i+1} = x^j y^{d-j+1}$, 
$b'_{j-i+2} = x^{j+1} y^{d-j}$, and this concludes the proof.
\end{proof}

Part~(b) of this lemma can be visualized as follows.

\begin{example}\label{ex-visual}
In~$\mathbb{T}^2$, consider the two order ideals $\OO$ and $\OO'$
with $\mu=18$ and $d=5$ given in Figure~\ref{fig-MaxDeg18}.

\begin{figure}[ht]
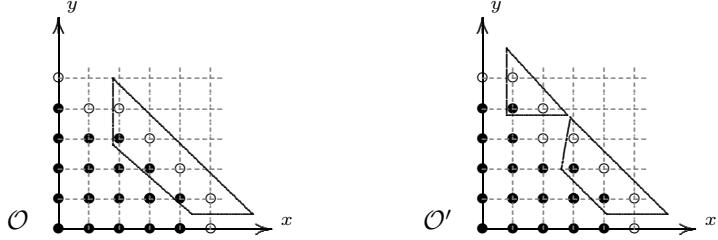

$\OO$ \
	\makebox{\beginpicture
		\setcoordinatesystem units <.4 cm, .4 cm>
		\setplotarea x from 0 to 7, y from 0 to 7
		\axis left /
		\axis bottom /
		\arrow <2mm> [.2,.67] from  6.5 0  to 7 0
		\arrow <2mm> [.2,.67] from  0 6.5  to 0 7
		\put {$\scriptstyle x$} [lt] <0.5mm,0.8mm> at 7.1 0.2
		\put {$\scriptstyle y$} [rb] <1.7mm,0.7mm> at 0.3 7
		\put {$\bullet$} at 0 0
		\put {$\bullet$} at 0 1	
		\put {$\bullet$} at 1 0
		\put {$\bullet$} at 1 1		
		\put {$\bullet$} at 2 0
		\put {$\bullet$} at 1 1
		\put {$\bullet$} at 0 2				
		\put {$\bullet$} at 3 0
		\put {$\bullet$} at 2 1
		\put {$\bullet$} at 1 2
		\put {$\bullet$} at 0 3
		\put {$\bullet$} at 0 4				
		\put {$\bullet$} at 4 0
		\put {$\bullet$} at 3 1
		\put {$\bullet$} at 1 3
		\put {$\bullet$} at 2 2
	 	\put {$\bullet$} at 2 3

		\put {$\bullet$} at 3 1
		\put {$\bullet$} at 4 1
		\put {$\bullet$} at 3 2
		
		\put {$\circ$} at 5 0
		\put {$\circ$} at 5 1		
		\put {$\circ$} at 3 1
		\put {$\circ$} at 4 0
		\put {$\circ$} at 0 5
		\put {$\circ$} at 1 4
		\put {$\circ$} at 2 4
		\put {$\circ$} at 3 3
		\put {$\circ$} at 4 2
\gray{	 		
		\setdashes<0.5mm>
		\plot 1 0  1 5.5 /
		\plot 2 0  2 5.5 /
		\plot 3 0  3 5.5 /
		\plot 4 0  4 5.5 /
		\plot 5 0  5 5.5 /
		\plot 0 1  5.5 1 /
		\plot 0 2  5.5 2 /
		\plot 0 3  5.5 3 /
		\plot 0 4  5.5 4 /
		\plot 0 5  5.5 5 /
}

\setlinear \plot 1.5  5.0    1.5 2.8    4.1  0.5    6.1  0.5    1.5  5.0 /

\endpicture} 
\qquad	\qquad $\OO'$ \	
\makebox{\beginpicture
		\setcoordinatesystem units <.4 cm, .4 cm>
		\setplotarea x from 0 to 7, y from 0 to 7
		\axis left /
		\axis bottom /
		\arrow <2mm> [.2,.67] from  6.5 0  to 7 0
		\arrow <2mm> [.2,.67] from  0 6.5  to 0 7
		\put {$\scriptstyle x$} [lt] <0.5mm,0.8mm> at 7.1 0.2
		\put {$\scriptstyle y$} [rb] <1.7mm,0.7mm> at 0.3 7
		
		\put {$\bullet$} at 0 0
		\put {$\bullet$} at 0 1	
		\put {$\bullet$} at 1 0
		\put {$\bullet$} at 1 1		
		\put {$\bullet$} at 2 0
		\put {$\bullet$} at 1 1
		\put {$\bullet$} at 0 2				
		\put {$\bullet$} at 3 0
		\put {$\bullet$} at 2 1
		\put {$\bullet$} at 1 2
		\put {$\bullet$} at 0 3
		
		\put {$\bullet$} at 0 4				
		\put {$\bullet$} at 4 0
		\put {$\bullet$} at 3 1
		\put {$\bullet$} at 1 3
		\put {$\bullet$} at 2 2
		\put {$\bullet$} at 3 1
		\put {$\bullet$} at 3 2
		\put {$\bullet$}  at 1 4
		\put {$\bullet$}  at 4 1
		
		\put {$\circ$} at 1 5
		\put {$\circ$} at 2 4
		\put {$\circ$} at 2 3
		\put {$\circ$} at 3 1
		\put {$\circ$} at 3 3
		\put {$\circ$} at 4 2
		\put {$\circ$} at 5 0
		\put {$\circ$} at 5 1		
  		\put {$\circ$} at 0 5
		
\gray{	 		
		\setdashes<0.5mm>
		\plot 1 0  1 5.5 /
		\plot 2 0  2 5.5 /
		\plot 3 0  3 5.5 /
		\plot 4 0  4 5.5 /
		\plot 5 0  5 5.5 /
		\plot 0 1  5.5 1 /
		\plot 0 2  5.5 2 /
		\plot 0 3  5.5 3 /
		\plot 0 4  5.5 4 /
		\plot 0 5  5.5 5 /
}

\setlinear \plot 0.5  6     0.5  3.8   2.5  3.8   0.5  6      /
\setlinear \plot 3.8  0.5   2.3  2.0   2.6  3.7   5.8  0.5  3.8  0.5  /

\endpicture}
\caption{Two order ideals with $d=5$ and their borders}\label{fig-MaxDeg18}
\end{figure}

Notice that $\OO$ has one $\OO_5$-segment, namely $(x^2y^3, x^3y^2, x^4y)$,
and~$\OO'$ has two $\OO_5$-segments, namely $(xy^4)$ and $(x^3y^2, x^4y)$. 

For the order ideal~$\OO$, we therefore get $\#\partial\OO = d+1+s =7$
and $\# C = \mu\nu = 18\cdot 7 = 126$. For the order ideal~$\OO'$, we find
$\#\partial\OO' = d+1+s = 8$ and $\#C = \mu\nu =  18 \cdot 8 = 144$.
\end{example}

This example indicates that, if we keep $\mu=\#\OO$ and $d=
\max\{\deg(t) \mid t\in\OO \}$ fixed, the number of indeterminates of~$K[C]$
increases with the segmentation type. Consequently, if we want to
construct a $Z$-separating re-embedding that shows $\BO \cong \mathbb{A}^{2\mu}$,
we have to eliminate more and more indeterminates. This leads us to the following
conjecture which we have verified for a large number of explicitly computed examples.

\begin{conjecture}\label{conj-newconj}
Let $\OO$ be a non-simplicial MaxDeg planar order ideal. Then the following conditions are equivalent.
\begin{enumerate}
\item There exists a $Z$-separating re-embedding of~$\BO$
which yields an isomorphism $\BO \cong \mathbb{A}^{2\mu}$.

\item  The segmentation type of~$\OO$ is one.
\end{enumerate} 
\end{conjecture}

Recall that all non-exposed indeterminates of~$C$ can be eliminated by Theorem~\ref{thm-elimNex}.
Sometimes this is already enough to construct the desired isomorphism.

\begin{example}\label{ex-122-again}
Let us reconsider the order ideal $\OO = \{1,y,x,y^2,xy\}$ given in Example~\ref{ex-122-2}.

\begin{figure}[ht]
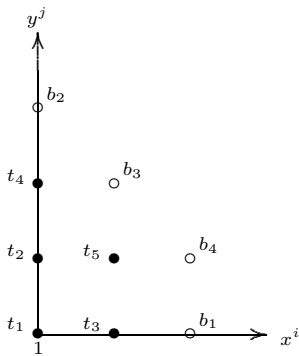

\centering{\makebox{
\beginpicture
		\setcoordinatesystem units <1cm,1cm>
		\setplotarea x from 0 to 3, y from 0 to 3.5
		\axis left /
		\axis bottom /
		\arrow <2mm> [.2,.67] from  2.5 0  to 3 0
		\arrow <2mm> [.2,.67] from  0 3.5  to 0 4
		\put {$\scriptstyle x^i$} [lt] <0.5mm,0.8mm> at 3.1 0
		\put {$\scriptstyle y^j$} [rb] <1.7mm,0.7mm> at 0 4
		\put {$\bullet$} at 0 0
		\put {$\bullet$} at 1 0
		\put {$\bullet$} at 0 1
		\put {$\bullet$} at 1 1
		\put {$\bullet$} at 0 2
		\put {$\scriptstyle 1$} [lt] <-1mm,-1mm> at 0 0
		\put {$\scriptstyle t_1$} [rb] <-1.3mm,0.4mm> at 0 0
		\put {$\scriptstyle t_3$} [rb] <-1.3mm,0.4mm> at 1 0
		\put {$\scriptstyle t_2$} [rb] <-1.3mm,0mm> at 0 1
		\put {$\scriptstyle t_4$} [rb] <-1.3mm,0mm> at 0 2
		\put {$\scriptstyle t_5$} [rb] <-1.3mm,0mm> at 1 1
		\put {$\scriptstyle b_1$} [lb] <0.8mm,0.8mm> at 2 0
		\put {$\scriptstyle b_2$} [lb] <0.8mm,0.8mm> at  0 3
		\put {$\scriptstyle b_3$} [lb] <0.8mm,0.8mm> at 1 2
	    \put {$\scriptstyle b_4$} [lb] <0.8mm,0.8mm> at 2 1
	
		\put {$\circ$} at 2 0
		\put {$\circ$} at 2 1
		\put {$\circ$} at 1 2
		\put {$\circ$} at 0 3
\endpicture
}} 
\caption{An order ideal with one segment}
\end{figure}

When the terms $t_3$, $t_4$, and~$t_5$ are multiplied by~$x$, the
result is in the border. The border terms $b_j$ such that $xb_j$ is involved
in a neighbor pair are~$b_2$ and~$b_3$. Hence the $x$-exposed indeterminates
are $c_{32}, c_{33}, c_{42}, c_{43}, c_{52}, c_{53}$. Similarly, we see
that the $y$-exposed indeterminates are $c_{41}, c_{43}, c_{44}, c_{51}, c_{53}, c_{54}$.
Altogether, there are 10 exposed indeterminates.

As we have $\# C = 20$ and $\dim(\BO)=10$, we want to eliminate 10 indeterminates.
By Theorem~\ref{thm-elimNex}, the tuple of non-exposed indeterminates
$$
Z = (c_{11}, c_{12}, c_{13}, c_{14}, c_{21}, c_{22}, c_{23}, c_{24}, c_{31}, c_{34})
$$
does the job, i.e., the $Z$-separating re-embedding yields an isomorphism
$\BO \cong \mathbb{A}^{10}$.
\end{example}

For an example in which the non-exposed indeterminates do not suffice, it is
enough to look at Example~\ref{ex-1232}. 

\begin{example}
The order ideal $\OO = \{1,y,x,y^2,xy,x^2,y^3,xy^2\}$ has a border basis scheme 
with an optimal $Z$-separating re-embedding $\BO \cong \mathbb{A}^{16}$, as shown 
in Example~\ref{ex-1232}. Since $\# C = 40$, it is not enough to eliminate
the 22 non-exposed indeterminates. Instead, one has to adjoin two exposed indeterminates
to~$Z$ to get the optimal re-embedding, as explained in Example~\ref{ex-1232}.
\end{example}

Finally, to see that planar MaxDeg order ideals with a segmentation type of two or more
have in general no optimal $Z$-separating re-embeddings, we have to look no further than
at the L-shape order ideal $\OO = \{1,y,x,y^2,x^2\}$. As we saw in Proposition~\ref{prop-L-Shape-Zsep},
the best possible $Z$-separating re-embeddings yield embeddings of~$\BO$ as codimension~2
subschemes of~$\mathbb{A}^{12}$, and the technique based on the Unimodular Matrix Problem is
required to construct an isomorphism with $\mathbb{A}^{10}$.

\bigskip

\subsection*{Acknowledgements.} The authors are grateful to L.N.\ Long for providing them with
the \cocoa\ package {\tt QuillenSuslin.cpkg5} which assisted the calculations for
Proposition~\ref{prop-theLshape}. The second authors thanks the University of Passau for
its hospitality during part of the preparation of this paper.

\bigskip\bigbreak
%
%


\begin{thebibliography}{99}

\bibitem{AKL} B.\ Andraschko, M.\ Kreuzer, and L.N.\ Long, Efficiently checking 
separating indeterminates, preprint 2024, 28 pages,
\url{https://doi.org/10.48550/arXiv.2412.18369}

\bibitem{Fog} J.\ Fogarty, Algebraic families on an algebraic surface,
Amer.\ J.\ Math.\ {\bf 90} (1968), 511--521, \url{https://doi.org/10.2307/2373541}

\bibitem{Gro} A.\ Grothendieck, Techniques de construction et th\'{e}or\`{e}mes 
d'existence en g\'{e}om\'{e}trie al\-g\'e\-bri\-que. IV. Les sch\'emas de Hilbert, 
Seminaire Bourbaki {\bf 221} (1961). Reprinted in: A.\ Douady et al. (eds), 
Seminaire Bourbaki, Vol. 6, Soc.\ Math.\ France, Paris 1995, pp.\ 249--276.

\bibitem{Hai} M.\ Haiman, $t,q$-Catalan numbers and the Hilbert schemes,
Discrete Math.\ {\bf 193} (1998), 201--224,
\url{https://doi.org/10.1016/S0012-365X%2898%2900141-1}

\bibitem{Hui1} M.\ Huibregtse, A description of certain affine open schemes
that form an open covering of ${\rm Hilb}_{\AA_2^k}^n$,
Pacific J.\ Math.\ {\bf 204} (2002), 97--143,
\url{https://doi.org/10.2140/pjm.2002.204.97}

\bibitem{Hui2} M.\ Huibregtse, The cotangent space at a monomial ideal
of the Hilbert scheme of points of an affine space, preprint 2005, 51 pages, 	
\url{https://doi.org/10.48550/arXiv.math/0506575}

\bibitem{Hui3} M.\ Huibregtse, An elementary construction of the multigraded 
Hilbert scheme of points, Pacific J.\ Math.\ {\bf 223} (2006), 269--315,
\url{https://doi.org/10.2140/pjm.2006.223.269}


\bibitem{KLR0} M.\ Kreuzer, L.N.\ Long, and L.\ Robbiano, 
Computing subschemes of the border basis scheme, 
Int. J. Algebra Comput. {\bf 30} (2020), 1671--1716,
\url{https://doi.org/10.1142/S0218196720500599}

\bibitem{KLR1} M.\ Kreuzer, L.N.\ Long, and L.\ Robbiano, 
Cotangent spaces and separating re-embeddings, J.\ Algebra Appl. {\bf 21} (2022), 
paper 2250188, \url{https://doi.org/10.1142/S0219498822501882}

\bibitem{KLR2} M.\ Kreuzer, L.N.\ Long, and L.\ Robbiano, 
Restricted Gr\"obner fans and re-embeddings of affine algebras, 
S\~ao Paulo J.\ Math.\ Sci.\ {\bf 17} (2023), 242--267, 
\url{https://doi.org/10.1007/s40863-022-00324-w}

\bibitem{KLR3} M.\ Kreuzer, L.N.\ Long, and L.\ Robbiano, 
Re-embeddings of affine algebras via Gr{\"o}bner fans of linear ideals, 
Beitr.\ Algebra Geom.\ {\bf 65} (2024), 827--851,
\url{https://doi.org/10.1007/s13366-024-00733-2}


\bibitem{KR1} M.\ Kreuzer and L.\ Robbiano, {\it Computational
Commutative  Algebra 1}, Springer-Verlag, Berlin Heidelberg, 2000,
\url{https://doi.org/10.1007/978-3-540-70628-1}

\bibitem{KR2} M.\ Kreuzer and L.\ Robbiano, {\it Computational
Commutative Algebra 2}, Springer-Verlag, Berlin Heidelberg, 2005,
\url{https://doi.org/10.1007/3-540-28296-3}

\bibitem{KR3} M.\ Kreuzer and L.\ Robbiano, Deformations of border
bases, Collect.\ Math.\ {\bf 59} (2008), 275--297,
\url{https://doi.org/10.1007/BF03191188}

\bibitem{KR4} M.\ Kreuzer and L.\ Robbiano, The geometry of border
bases, J.\ Pure Appl.\ Algebra {\bf 215} (2011), 2005--2018,
\url{https://doi.org/10.1016/j.jpaa.2010.11.011}

\bibitem{KR5} M.\ Kreuzer and L.\ Robbiano, Elimination by substitution,
preprint 2024, 26 pages, \url{https://doi.org/10.48550/arXiv.2403.06415}

\bibitem{KSL} M. Kreuzer, B.\ Sipal, and L.N.\ Long,
On the regularity of the monomial point of a border basis scheme,
Beitr. Algebra Geom. {\bf 61} (2020), 515--532,
\url{https://doi.org/10.1007/s13366-019-00482-7}




\end{thebibliography}
\end{document}